\documentclass[11pt]{amsart}
\usepackage{amsmath}
\usepackage{amsfonts}
\usepackage{amssymb}
\usepackage[arrow,matrix,curve,cmtip,ps]{xy}

\numberwithin{equation}{section}

\usepackage{amsthm}

\newtheorem{theorem}{Theorem}[section]
\newtheorem{lemma}[theorem]{Lemma}
\newtheorem{proposition}[theorem]{Proposition}
\newtheorem{corollary}[theorem]{Corollary}
\newtheorem{definition}[theorem]{Definition}

\newtheorem{remark}[theorem]{Remark}

\newcommand\id{\mathop{\rm id}}

\newcommand\nph{\varphi}

\newcommand\cp{\mathop{\rm cp}}

\newcommand\ttt{\mathop{\rm t}}

\newcommand\linj{\mathop{\rm el}}
\newcommand\rinj{\mathop{\rm er}}
\newcommand\inj{\mathop{\rm e}}
\newcommand\comm{\mathop{\rm c}}

\newcommand\OMIN{\mathop{\rm OMIN}}
 \newcommand\OMAX{\mathop{\rm OMAX}}

\newcommand{\cl}[1]{\mathcal{#1}}
\newcommand{\bb}[1]{\mathbb{#1}}

\begin{document}

\title{Tensor products of operator systems}

\author[A.~Kavruk]{Ali Kavruk}
\address{Department of Mathematics, University of Houston,
Houston, Texas 77204-3476, U.S.A.}
\email{kavruk@math.uh.edu}

\author[V.~I.~Paulsen]{Vern I.~Paulsen}
\address{Department of Mathematics, University of Houston,
Houston, Texas 77204-3476, U.S.A.}
\email{vern@math.uh.edu}

\author[I.~G.~Todorov]{Ivan G.~Todorov}
\address{Department of Pure Mathematics, Queen's University Belfast,
Belfast BT7 1NN, United Kingdom}
\email{i.todorov@qub.ac.uk}

\author[M.~Tomforde]{Mark Tomforde}
\address{Department of Mathematics, University of Houston, Houston,
  Texas 77204-3476, U.S.A.}
\email{tomforde@math.uh.edu}

\date{December 12, 2009, revised February 10, 2011}
\thanks{The first and second authors were supported by NSF grant
  DMS-0600191. The third author was supported by EPSRC grant D050677/1.
  The fourth author was supported by NSA Grant H98230-09-1-0036.}

\begin{abstract}
The purpose of the present paper is to lay the foundations for a
systematic study of tensor products of operator systems.
After giving an axiomatic definition of tensor products in this category, we examine in detail
several particular examples of tensor products, including a minimal,
maximal, maximal commuting, maximal injective and some asymmetric
tensor products. We characterize these tensor products in terms of
their universal properties and give descriptions of their positive cones.
We also characterize the corresponding tensor products of operator
spaces induced by a certain canonical inclusion of an operator space
into an operator system.
We examine notions of nuclearity for our tensor products which, on the category of C*-algebras, reduce to
the classical notion. We exhibit an operator system $\cl S$ which is not
completely order isomorphic to a C*-algebra yet has the property that
for every C*-algebra $A$,
the minimal and maximal tensor product of $\cl S$ and $A$ are equal.
\end{abstract}

\maketitle

\section{Introduction}

For the last 25 years there has been a great deal of development of
the theory of tensor products of operator spaces and there has been
a great influx of ideas and techniques from Banach space theory.
During the same period there has been very little development of the
tensor theory of operator systems. Since the methods of \cite{Pa}
show that many of the basic results about operator spaces and
completely bounded maps can be derived from results about operator
systems and completely positive maps, we believe that further
development of the tensor theory of operator systems could play an
important role in operator space tensor theory, as well as having
its own intrinsic merit.

In this paper we initiate the systematic study of tensor products
in the category whose objects are operator systems and whose morphisms are
unital completely positive maps.
After setting the axiomatic foundations in Section 3,
we introduce and study several particular tensor products.
We thus dedicate Section 4 to the \lq\lq minimal'' tensor product of operator systems,
which corresponds to the formation of composite
quantum systems in Quantum Information Theory.

In Sections 5 and 6 we study
the \lq\lq maximal'' and the \lq\lq commuting'' tensor products.
The maximal tensor product is also important in Quantum Information
Theory, since the states
on the minimal tensor product of two finite dimensional operator
systems can be identified with the maximal tensor product of their
dual spaces \cite[Proposition~1.9]{FP}.
We characterize the maximal
tensor product in terms of a universal linearization property for jointly
completely positive maps, and the commuting tensor product in terms
of the maximal C*-algebraic tensor product of certain universal
C*-algebras associated with the corresponding operator systems. It
follows from an earlier work of Lance \cite{La} that, given two
C*-algebras, their maximal and commuting tensor products as operator systems both agree with their
C*-maximal tensor product. However, we show that for general
operator systems these tensor products are distinct. Thus the
maximal tensor product and the commuting tensor product give two
different ways to extend the C*-maximal tensor product from the
category of C*-algebras to operator systems. This implies that
C*-algebraic notions that can be defined in terms of the minimal and
maximal C*-tensor products, such as nuclearity, weak expectation
property (WEP), and exactness, can bifurcate into multiple concepts
in this larger category.

In particular, we exhibit an operator system $\cl S$ which is not
completely order isomorphic to a C*-algebra and which does not
``factor through matrix algebras''; i.e., is not nuclear in this classical sense, but which has
the property that for every C*-algebra $A$, the minimal and the
maximal operator system tensor product structures on $\cl S \otimes
A$ coincide. Similarly, we exhibit operator systems that are not
nuclear in the classical sense, but which have the property that
their minimal and commuting tensor products with every operator
system are equal. This is achieved through a careful examination of operator subsystems
of the space of all $n$ by $n$ matrices associated with graphs.

Since every operator space embeds in a canonical operator
system, tensor products in the operator system category can be
pulled back to tensor products in the operator space category. We
describe the pullbacks of the operator system tensor products that
we construct. In particular, we show that the tensor product induced
by the maximal (respectively,
minimal) operator system tensor product coincides with the operator
projective (respectively, injective) tensor product. The family of
tensor products on the operator space category that one can obtain
as pullbacks is potentially more suited for carrying out
Grothendieck's program.

In Section 7, we examine the lattice structure of operator system tensor products.
This allows us to introduce maximal, one-sided and two-sided, injective tensor products,
the one-sided ones being asymmetric. We also formulate a characterization of nuclearity and WEP for
C*-algebras in terms of these asymmetric tensor products.

\section{Preliminaries}

In this section we establish the terminology and state the definitions that shall be used
throughout the paper.

A {\bf $*$-vector space} is a complex vector space $V$ together with
a map $^* : V \to V$ that is involutive (i.e., $(v^*)^* = v$ for
all $v \in V$) and conjugate linear (i.e., $(\lambda v + w)^* =
\overline{\lambda} v^* + w^*$ for all $\lambda \in \bb C$ and $v,w
\in V$). If $V$ is a $*$-vector space, then we let $V_h = \{x \in V
: x^* = x \}$ and we call the elements of $V_h$ the {\bf hermitian}
elements of $V$. Note that $V_h$ is a real vector space.

An {\bf ordered $*$-vector space} is a pair $(V, V^+)$ consisting of a $*$-vector space $V$ and a subset $V^+ \subseteq V_h$ satisfying
the following two properties:

\medskip

(a) $V^+$ is a cone in $V_h$;

(b) $V^+ \cap -V^+ = \{ 0 \}$.

\medskip

In any ordered $*$-vector space we may define a partial order $\geq$ on $V_h$
by defining $v \geq w$ (or, equivalently, $w \leq v$) if and only if $v - w \in V^+$.
Note that $v \in V^+$ if and only if $v
\geq 0.$ For this reason $V^+$ is called the cone of
{\bf positive} elements of $V$.

If $(V,V^+)$ is an ordered $*$-vector space, an element
$e \in V_h$ is called an {\bf order unit} for $V$ if
for all $v \in V_h$ there exists a real number $r > 0$
such that $re \geq v$.
If $(V, V^+)$ is an ordered $*$-vector space with an order unit $e$,
then we say that $e$ is an {\bf Archimedean order unit} if whenever
$v \in V$ and $re+v \geq 0$ for all real $r >0$, we have that $v \in V^+$.
In this case, we call the triple $(V, V^+, e)$ an {\bf Archimedean
ordered $*$-vector space} or an {\bf AOU space,} for short.
The {\bf state space} of $V$ is the set $S(V)$ of all linear maps $f : V\rightarrow\bb{C}$
such that $f(V^+)\subseteq [0,\infty)$ and $f(e) = 1$.

If $V$ is a $*$-vector space, we let $M_{m,n}(V)$ denote the set of
all $m \times n$ matrices with entries in $V$ and set $M_{n}(V) =
M_{n,n}(V)$.  The natural
addition and scalar multiplication turn $M_{m,n}(V)$ into a complex
vector space. We set $M_{m,n} := M_{m,n}(\bb{C})$, and let
$\{E_{i,j} : 1\leq i\leq n, 1\leq j \leq m\}$ denote its canonical
matrix unit system. If $X = (x_{i,j})_{i,j} \in M_{l,m}$ is a scalar
matrix, then for any $A = (a_{i,j})_{i,j} \in M_{m,n}(V)$ we let
$XA$ be the element of $M_{l,n}(V)$ whose $i,j$-entry $(XA)_{i,j}$
equals $\sum_{k=1}^m x_{i,k} a_{k,j}$. We define multiplication by
scalar matrices on the left in a similar way. Furthermore, when
$m=n$, we define a $*$-operation on $M_n(V)$ by letting
$(a_{i,j})_{i,j}^* := (a_{j,i}^*)_{i,j}$. With respect to this
operation, $M_n(V)$ is a $*$-vector space. We let $M_n(V)_h$ be the
set of all hermitian elements of $M_n(V)$.

\begin{definition}
Let $V$ be a $*$-vector space.  We say that $\{ C_n \}_{n=1}^\infty$ is a \emph{matrix ordering} on $V$ if
\begin{enumerate}
\item $C_n$ is a cone in $M_n(V)_h$ for each $n\in \bb{N}$,
\item $C_n \cap -C_n = \{0 \}$ for each $n\in \bb{N}$, and
\item for each $n,m\in \bb{N}$ and $X \in M_{n,m}$
we have that $X^* C_n X \subseteq C_m$.
\end{enumerate}
In this case we call $(V, \{C_n \}_{n=1}^\infty )$ a \emph{ matrix
  ordered $*$-vector space.} We refer to condition (3) as the
  \emph{compatibility of the family
  $\{C_n\}_{n=1}^{\infty}$}.
\end{definition}

Note that properties (1) and (2) show that $(M_n(V), C_n)$
is an ordered $*$-vector space for each $n\in \bb{N}$. As usual, when $A,B \in M_n(V)_h$, we
write $A \leq B$ if $B-A \in C_n$.

\begin{definition}\label{op-system-abstract-def}
Let $(V, \{ C_n \}_{n=1}^\infty)$ be a matrix ordered $*$-vector space.  For $e \in V_h$ let
$$e_n := \left( \begin{smallmatrix} e & & \\ & \ddots & \\ & & e \end{smallmatrix} \right)$$
be the corresponding diagonal matrix in $M_n(V)$.
We say that $e$ is a \emph{matrix order unit} for $V$
if $e_n$ is an order unit for $(M_n(V), C_n)$ for each $n$.
We say that $e$ is an \emph{Archimedean matrix order unit} if
$e_n$ is an Archimedean order unit for $(M_n(V), C_n)$ for each $n$.
An \emph{(abstract) operator system} is a triple $(V, \{C_n \}_{n=1}^\infty, e)$,
where $V$ is a complex $*$-vector space,
$\{C_n \}_{n=1}^\infty$ is a matrix ordering on $V$,
and $e \in  V_h$ is an Archimedean matrix order unit.
\end{definition}

We note that the above definition of an operator system was first
introduced by Choi and Effros in \cite{CE2}.  If $V $ and $V'$ are
vector spaces and $\phi : V \to V'$ is a linear map, then for each
$n\in \bb{N}$ the map $\phi$ induces a linear map $\phi^{(n)} : M_n(V)
\to M_n(V')$ given by $\phi^{(n)} ((v_{i,j})_{i,j}) :=
(\phi(v_{i,j}))_{i,j}$. If $(V, \{ C_n \}_{n=1}^\infty)$ and $(V',
\{ C_n' \}_{n=1}^\infty)$ are matrix ordered $*$-vector spaces, a
map $\phi : V \to V'$ is called {\bf completely positive} (for short, {\bf c.p.}) if
$\phi^{(n)}(C_n) \subseteq C_n'$ for each $n\in \bb{N}$. Similarly, we
call a linear map $\phi : V \to V'$ a {\bf complete order
isomorphism} if $\phi$ is invertible and both $\phi$ and $\phi^{-1}$
are completely positive.

We denote by $\cl B(H)$ the space of all bounded linear operators
acting on a Hilbert space $H$.
The direct sum of $n$ copies of $H$ is denoted by $H^n$ and its elements are written as column vectors.
A {\bf concrete operator system}
$\cl S$ is a subspace of $\cl B (H)$ such that $\cl S = \cl S^*$ and
$I \in \cl S$. (Here, and in the sequel, we let $I$ denote the
identity operator.) As is the case for many classes of subspaces
(and subalgebras) of $\cl B (H)$, there is an abstract
characterization of concrete operator systems.
If $\cl S \subseteq \cl B (H)$
is a concrete operator system, then we observe that $\cl S$ is a
$*$-vector space with respect to the adjoint operation, $\cl S$ inherits an order structure from $\cl B(H)$,
and has $I$ as an Archimedean order unit. Moreover, since $\cl
S \subseteq \cl B(H)$, we have that $M_n(\cl S) \subseteq M_n( \cl B(H)) \equiv \cl B (H^n)$ and hence $M_n(\cl S)$ inherits
an involution and an order
structure from $\cl B (\cl H^n)$ and has the $n \times n$ diagonal matrix
$$\begin{pmatrix} I & & \\ & \ddots & \\ & & I \end{pmatrix}$$ as an
Archimedean order unit.  In other words, $\cl S$ is
an abstract operator system in the sense of Definition
\ref{op-system-abstract-def}. The following result of Choi and
Effros \cite[Theorem~4.4]{CE2} shows that the converse is also true.
For an alternative proof of the result, we refer the reader to
\cite[Theorem~13.1]{Pa}.

\begin{theorem}[Choi-Effros]
\label{th-choieffros}
Every concrete operator system $\cl S$ is an abstract operator system.
Conversely, if $(V, \{C_n \}_{n=1}^\infty, e)$ is an abstract operator system, then there
exists a Hilbert space $\cl H$, a concrete operator system
$\cl S \subseteq \cl B(H)$, and a complete order
isomorphism $\phi : V \to \cl S$ with $\phi(e) = I$.
\end{theorem}

Thanks to the above theorem, we can identify abstract and concrete
operator systems and refer to them simply as operator systems.
To avoid excessive notation, we will generally refer to an operator
system as simply a set $\cl S$ with the understanding that $e$ is the
order unit and $M_n(\cl S)^+$ is the cone of positive elements
in $M_n(\cl S)$. We note that any unital C*-algebra (and all C*-algebras in the
present paper will be assumed to be unital) is also an operator system in a canonical way.

There is a similar theory for arbitrary subspaces $X \subseteq \cl
B(H)$, called also {\bf concrete operator spaces}. The
identification $M_n(\cl B(H)) \equiv \cl B(H^n)$ endows each $M_n(X)
\subseteq M_n(\cl B(H))$ with a norm; the family of norms obtained
in this way satisfies certain compatibility axioms called {\it
Ruan's axioms.} Ruan's theorem identifies the vector spaces
satisfying Ruan's axioms with the concrete operator spaces. Sources
for the details include \cite{er2} and \cite{Pa}.

What is important for our setting is that the dual of every operator
space is again an operator space \cite{bp, er2} and that the dual of an operator
system is a matrix-ordered space \cite{CE2}. Thus the dual of an operator system
carries two structures and we will need to understand the
relationship between these structures.

To this end, let $\cl S$ be an operator system and let $\cl S^d$
denote its Banach space dual. For $f \in \cl S^d,$ we define $f^*
\in \cl S^d$ by $f^*(s) = \overline{f(s^*)}.$ This operation turns
$\cl S^d$ into a $*$-vector space and it is easy to check that the
cone of positive linear functionals defines an order on $\cl S^d.$
One can define a matrix order by declaring an element $(f_{i,j}) \in
M_n(\cl S^d)$ to be {\bf positive} if and only if the map $F : \cl S
\to M_n$ given by $F(s) = (f_{i,j}(s))$ is completely positive. It
follows from \cite[Lemma~4.2, Lemma 4.3]{CE2} that this family of
sets is a matrix ordering on $\cl S^d.$

On the other hand, one defines a norm on $M_n(\cl S^d)$ by setting
$\|(f_{i,j})\| = \|F\|_{cb},$ where $\|F\|_{cb}$ denotes the
completely bounded norm of the mapping $F.$  This family of norms
satisfies Ruan's axioms and thus
gives $\cl S^d$ the structure of an abstract operator space.

The following result compares these two structures.

\begin{theorem} Let $\cl S$ be an operator system. Then there exists a
  Hilbert space $H$ and a weak* continuous completely positive map $\Phi: \cl S^d \to
  \cl B(H)$ that is a complete order isomorphism onto its range and
  satisfies
$$\|(\Phi(f_{i,j}))\| \le \|(f_{i,j})\| \le 2 \|(\Phi(f_{i,j}))\|$$
for all $(f_{i,j}) \in M_n(\cl S^d)$ and all $n\in \bb{N}$.
\end{theorem}
\begin{proof} Let $\cl I_n = \{P \in M_n(\cl S)^+: \|P\|
  \le 1 \},$ so that $0 \le P \le e_n$ for each $P\in \cl I_n$.
  For each $P = (p_{i,j})\in \cl I_n$ define
  $e_P: \cl S^d \to M_n$ by setting $e_P(f) = (f(p_{i,j})).$ The
  map $e_P$ is completely positive by \cite[Lemma 4.3]{CE2} and since $\|P\| \le 1,$ we have
  that $\|e_P\|_{cb} \le 1.$ Note that the space $A_n = \ell^{\infty}(\cl I_n, M_n)$
  of all bounded $M_n$-valued functions defined on the set $\cl I_n$
  is a unital C*-algebra and that $M_k(A_n) \equiv \ell^{\infty}(\cl I_n,M_{kn})$ in a canonical way.
  Let $\phi_n: \cl S^d \to A_n$ be defined by
  $\phi_n(f)(P) = e_P(f)$. It follows that $\phi_n$ is completely positive and $\|\phi_n\|_{cb}
  \le 1.$

Now define $\Phi: \cl S^d \to \sum_{n=1}^{\infty} \oplus A_n$ by letting $\Phi(f)
= \sum_{n=1}^{\infty} \oplus \phi_n(f)$; we have that $\Phi$ is completely positive
and $\|\Phi\|_{cb} \le 1$. Since $(f_{i,j}) \in M_n(\cl
S^d)^+$ if and only if $(e_P(f_{i,j})) \ge 0$ for every $P \in \cl
I_m$ and every $m,$ we have that $\Phi$ is a complete order
isomorphism onto its range.  It is also clear that $\Phi$ is weak* continuous.

Let $(f_{i,j})\in M_n(\cl S^d)$ and $F : \cl S\rightarrow M_n$ be the map given by
$F(s) = (f_{i,j}(s))$.
Given any $x \in M_n(\cl S)$ with $\|x\| \le 1$
we have that
\[ P = \frac{1}{2} \begin{pmatrix} e_n & x\\ x^* &
  e_n \end{pmatrix} \in \cl I_{2n}, \]
and hence $\frac{1}{2}\|F(x)\| \le \|(e_P(f_{i,j}))\| \le
\|(\Phi(f_{i,j}))\|.$ Thus, $\|(f_{i,j})\| = \|F\|_{cb} \le 2
\|(\Phi(f_{i,j}))\|,$ and the result follows.
\end{proof}


Given two operator systems, $\cl S$ and $\cl T$, we write
$\operatorname{CP} (\cl S,\cl T)$ for the cone of all completely
positive maps from $\cl S$ into $\cl T$, and we write
$\operatorname{UCP} (\cl S,\cl T)$ for the set of all unital
completely positive (abbreviated u.c.p.) maps from $\cl S$ into $\cl
T$. We denote by {\bf $\cl{O}$} the category whose objects are
operator systems and whose morphisms are unital completely positive
maps. The {\bf matricial state space} of an operator system $\cl S$
is the set $S_{\infty}(\cl S) = \bigcup_{n=1}^{\infty} S_n(\cl S)$,
where
$$S_n(\cl S) = \{\phi : \cl S\rightarrow M_n \ : \ \phi \mbox{ a
unital completely positive map}\}.$$

The algebraic tensor product of two vector spaces $V$ and $W$ is
denoted by $V\otimes W$. If $V^+\subseteq V$ and $W^+\subseteq W$
are cones, we let $V^+\otimes W^+ = \{v\otimes w : v\in V^+, w\in
W^+\}$. For $n,m\in \bb{N}$, we shall use the usual Kronecker
identification of $M_n \otimes M_m$ with $ M_{mn}$; thus, if
$(x_{i,j}) \in M_n$ and $(y_{k,l}) \in M_m$, we identify $(x_{i,j})
\otimes (y_{k,l})$ with the matrix
$(x_{i,j}y_{k,l})_{(i,k),(j,l)}\in M_{mn}.$ At the level of matrix
units we have $E_{i,j} \otimes E_{k,l} = E_{(i,k),(j,l)}.$


If $V_1,V_2$, and $W$ are vector spaces and if $\psi : V_1\times
V_2\rightarrow W$ is a bilinear map, then for $n,m\in \bb{N}$ we let
$\psi^{(n,m)} : M_n(V_1)\times M_m(V_2)\rightarrow M_{nm}(W)$ be the
bilinear map
given by $\psi^{(n,m)}((x_{i,j})_{i,j},(y_{k,l})_{k,l}) = (\psi(x_{i,j},y_{k,l}))_{(i,k),(j,l)}$.

Another construction that will play a role throughout this paper is
the {\bf Archimedeanization} of an ordered (respectively, matrix
ordered) $*$-vector space with an order unit $e.$ This was first
introduced in \cite{pt} for ordered spaces and extended to matrix
ordered spaces in \cite{ptt}.  Briefly, if $(V,
\{D_n\}_{n=1}^{\infty},e)$ is a matrix ordered $*$-vector space with
matrix order unit $e$ and with the property that $(V, D_1, e)$ is an
AOU space, then the Archimedeanization is obtained by forming the
smallest set of cones $C_n \subseteq M_n(V),$ such that
$D_n\subseteq C_n$ and $(V, \{C_n \}_{n=1}^{\infty}, e)$ is an
operator system.  In \cite{ptt}, an explicit description of the
elements of $C_n$ is given; namely, we have that $C_n = \{p\in
M_n(V) : p + re_n \in D_n, \mbox{ for all } r > 0\}$.
 We record one fact
about this process that we shall need later.

\begin{lemma}\label{l-cparch}
Let $(V, \{D_n\}_{n=1}^{\infty},e)$ be a matrix ordered
$*$-vector space with matrix order unit $e$ and with the property that $(V, D_1, e)$ is an AOU space. Let $(C_n)_{n=1}^{\infty}$ be the cones obtained through the Archimedeanization process.
Suppose that $\cl T$ is an operator system and $\phi : V\rightarrow \cl T$ is a linear map. We have that $\phi^{(n)}(D_n)\in M_n(\cl T)^+$
if and only if $\phi^{(n)}(C_n)\in M_n(\cl T)^+$, for each $n\in \bb{N}$.
\end{lemma}
\begin{proof}
This follows from the characterization of the Archmedeanization as the
smallest set of cones turning $V$ into an operator system.
\end{proof}

We shall also frequently need the following fact.

\begin{lemma}\label{l-detoss}
Let $V$ be a vector space and $\cl S$ and $\cl T$ be operator systems
with underlying vector space $V$. Suppose that
$\operatorname{UCP} (\cl S,$ $\cl B(H)) = \operatorname{UCP} (\cl T,\cl B(H))$ for every Hilbert space $H$.
Then $\cl S$ is completely order isomorphic to $\cl T$.
\end{lemma}
\proof
Assume, without loss of generality, that $\cl S\subseteq \cl B(H)$ is a concrete operator system. Then
the identity map $\id: \cl S\rightarrow\cl B(H)$ is unital and completely positive.
It follows that $\id$ is completely positive on $\cl T$ and hence $M_n(\cl T)^+\subseteq M_n(\cl S)^+$.
Reversing the argument implies that the identity map on $V$ is a unital complete
order isomorphism.
\endproof


\section{Tensor Products of Operator Systems}\label{s-gs}

We start this section with the definitions of the main concepts
studied in this paper.
Given a pair of operator systems
$(\cl S, \{ P_n\}_{n=1}^{\infty}, e_1)$ and $ (\cl T,\{
Q_n\}_{n=1}^{\infty}, e_2)$ by an {\bf operator system structure
on $\cl S \otimes \cl T$,} we mean a family $\tau = \{C_n\}_{n=1}^{\infty}$
of cones, where $C_n \subseteq M_n(\cl S \otimes \cl T)$, satisfying:

\begin{itemize}
\item[(T1)] $(\cl S \otimes \cl T,\{C_n\}_{n=1}^{\infty}, e_1 \otimes e_2)$ is
an operator system denoted $\cl S \otimes_{\tau} \cl T$,
\item[(T2)] $P_n \otimes Q_m \subseteq C_{nm},$ for all $n,m \in \bb
N$, and
\item[(T3)] If $\phi:\cl S \to M_n$  and $\psi:\cl T \to M_m$ are unital completely
positive maps, then $\phi \otimes \psi:\cl S \otimes_{\tau} \cl T \to M_{mn}$ is a unital completely
positive map.
\end{itemize}

To simplify notation we shall generally write $C_n = M_n(\cl S \otimes_{\tau} \cl T)^+.$
Conditions (T2) and (T3) are reminiscents of Grothendieck's axioms for tensor products of Banach spaces.
Condition (T2) may be viewed as the order analogue of the cross-norm condition,
while (T3) as the analogue of the property of a cross-norm of being ``reasonable''.

Given two operator system structures $\tau_1$
and $\tau_2$ on $\cl S \otimes\cl T,$ we say that
{\bf $\tau_1$ is greater than $\tau_2$} provided that the identity map
on $\cl S \otimes \cl T$ is completely positive from $\cl S
\otimes_{\tau_1} \cl T$ to $\cl S \otimes_{\tau_2} \cl T,$ which is
equivalent to requiring that $M_n(\cl S \otimes_{\tau_1} \cl T)^+
\subseteq M_n(\cl S \otimes_{\tau_2} \cl T)^+$ for every $n \in \bb N.$

By an {\bf operator system tensor product,} we mean a
mapping  $\tau: \cl O \times \cl O \to \cl O,$
such that for every pair of operator systems $\cl S$ and $\cl T,$
 $\tau(\cl S, \cl T)$ is an operator system structure on $\cl
S \otimes \cl T,$ denoted $\cl S\otimes_{\tau}\cl T$.

We call an operator system tensor product $\tau$ {\bf functorial,}
if the following property is satisfied:

\begin{itemize}
\item[(T4)] For any four operator systems $\cl S_1,\cl S_2, \cl T_1$, and $\cl T_2$, we have that if $\phi\in \operatorname{UCP} (\cl S_1,\cl S_2)$ and $\psi\in \operatorname{UCP}(\cl T_1,\cl
T_2)$, then the linear map $\phi\otimes\psi : \cl S_1\otimes\cl
T_1\rightarrow\cl S_2\otimes\cl T_2$ belongs to
$\operatorname{UCP}(\cl S_1\otimes_{\tau}\cl T_1,\cl
S_2\otimes_{\tau}\cl T_2)$.
\end{itemize}

If for all operator systems $\cl S$ and $\cl T$ the map $\theta :
x\otimes y\rightarrow y\otimes x$ extends to a unital complete order
isomorphism from $\cl S\otimes_{\tau}\cl T$ onto $\cl
T\otimes_{\tau}\cl S$ then $\tau$ is called {\bf symmetric}.

Given three vector spaces $\cl R, \cl S,$ and $\cl T$, there is a
natural isomorphism from $(\cl R \otimes \cl S) \otimes \cl T$ onto
$\cl R \otimes ( \cl S \otimes \cl T).$ We say that an operator system
tensor product $\tau$ is {\bf associative} if for any three
operator systems $\cl R, \cl S,$ and $\cl T$, this natural isomorphism
yields a complete order isomorphism from $(\cl R\otimes_{\tau} \cl S)\otimes_{\tau}
\cl T$ onto $\cl R\otimes_{\tau}(\cl S\otimes_{\tau}\cl T).$

We say that a functorial operator system tensor product is {\bf injective} if
for all operator systems $\cl S_1\subseteq\cl S_2$ and $\cl
T_1\subseteq\cl T_2$, the inclusion $\cl S_1\otimes_{\tau}\cl T_1
\subseteq \cl S_2\otimes_{\tau}\cl T_2$ is a complete order
isomorphism onto its range, that is,
$M_n(\cl S_1\otimes \cl T_1)\cap M_n(\cl S_2\otimes_{\tau} \cl T_2)^+ = M_n(\cl S_1\otimes_{\tau} \cl T_1)^+$
for every $n\in \bb{N}$.

One important concept from the theory of C*-algebras that we shall be
interested in generalizing is {\it nuclearity}.

\begin{definition}\label{d_alpbet}
Let $\alpha$ and $\beta$ be operator system tensor
products. An operator system $\cl S$ will be called
$(\alpha,\beta)$-nuclear if the identity map between $\cl
S\otimes_{\alpha}\cl T$ and $\cl
S\otimes_{\beta}\cl T$ is a complete order isomorphism for every operator system $\cl T$.
\end{definition}

One shortcoming of the theory of operator space tensor products is
that the minimal and maximal operator space tensor products of matrix
algebras do not coincide. For this reason there are essentially no
nuclear spaces in the operator space category. We will see that, unlike the operator space
case, there is a rich theory of nuclear operator systems for the
various tensor products we will introduce subsequently.

Recall that every operator system is also an operator space whose
matrix norms are determined by the matrix order.  In fact, if $\cl
S$ is an operator system with order unit $e,$ then $s=(s_{i,j}) \in
M_n(\cl S)$ satisfies $\|(s_{i,j})\| \le 1$ if and only if
$\begin{pmatrix} e_n & s\\s^* & e_n \end{pmatrix} \in M_2(M_n(\cl
S))^+.$ Since we shall need this fact often, it is worthwhile to
write it out in tensor notation.  Thus, we have that
$\|\sum_{i,j=1}^n E_{i,j} \otimes s_{i,j} \| \le 1$ if and only if
$E_{1,1} \otimes e_n +E_{2,2} \otimes e_n + E_{1,2} \otimes s +
E_{2,1} \otimes s^*= \sum_{i=1}^n (E_{1,1}+E_{2,2}) \otimes E_{i,i}
\otimes e + \sum_{i,j=1}^n (E_{1,2} \otimes E_{i,j} \otimes s_{i,j}
+ E_{2,1} \otimes E_{i,j} \otimes s_{j,i}^*)$ is in $(M_2 \otimes
M_n \otimes \cl S)^+= M_{2n}(\cl S)^+.$

Since operator systems are also operator spaces, it is important to
understand the relationship between operator system tensor products
and operator space tensor products. But first, we record some two elementary facts that will be useful throughout.

\begin{proposition}\label{non-unitalcp}
Let $\cl S$ and $\cl T$ be operator systems and let
$\tau$ be an operator system structure on $\cl S
  \otimes \cl T.$ If $\phi: \cl S \to M_n$ and $\psi: \cl T \to M_m$ are completely positive, then $\phi \otimes \psi: \cl S \otimes_{\tau} \cl T \to M_{mn}$ is completely positive.
\end{proposition}
\begin{proof} By \cite[Exercise~6.2]{Pa}, there exist unital
  completely positive maps $\phi_1: \cl S \to M_n$ and $\psi_1: \cl T
  \to M_m$ and positive matrices $P \in M_n, Q \in M_m$ such that
  $\phi(x) = P\phi_1(x)P$ and $\psi(y) = Q\psi_1(y)Q.$  Hence, $\phi
  \otimes \psi(x \otimes y) = (P \otimes Q)( \phi_1 \otimes \psi_1 (x
  \otimes y))(P \otimes Q).$  By Property (T3), $\phi_1 \otimes \psi_1: \cl S
  \otimes_{\tau} \cl T \to M_{mn}$ is completely positive, and the result follows.
\end{proof}

The next fact is a trick that is sometimes used in the theory of ``decomposable'' maps.

\begin{proposition}\label{tildetrick}
Let $\cl S$ and $\cl T$ be operator systems and let $\gamma_{i,j}: \cl S \to \cl T,$ $1 \le i,j \le n$ be linear maps. Define $\Gamma: \cl S \to M_n(\cl T)$ by $\Gamma(x) = ( \gamma_{i,j}(x))$ and $\widetilde{\Gamma}: M_n(\cl S) \to M_n(\cl T)$ by $\widetilde{\Gamma}((x_{i,j})) = (\gamma_{i,j}(x_{i,j})).$ Then $\Gamma$ is completely positive if and only if $\widetilde{\Gamma}$ is completely positive.
\end{proposition}
\begin{proof} First assume that $\widetilde{\Gamma}$ is completely positive. Since the map $\delta: \cl S \to M_n(\cl S)$ defined by $\delta(x) = (x_{i,j})$ where $x_{i,j} =x$ for all $1 \le i,j \le n,$ is completely positive and $\Gamma(x) = \widetilde{\Gamma} \circ \delta(x),$ it follows that $\Gamma$ is completely positive.

Conversely, if $\Gamma$ is completely positive, then the map
$\Gamma^{(n)}:M_n(\cl S) \to M_n(M_n(\cl T)),$ is completely
positive. The map defined by compressing a matrix in $M_n(M_n(\cl
T))$ to a matrix in $M_n(\cl T),$ by letting the $(i,j)$-th entry of
the latter to be equal to the $(i,j)$-th entry of the $(i,j)$-th block of the former,
is completely positive, and the composition of
$\Gamma^{(n)}$ with this compression equals $\widetilde{\Gamma}$.
More precisely,
identifying $M_n(M_n)$ with $\cl B(\bb C^n \otimes \bb C^n)$ and letting $V: \bb C^n
\to \bb C^n \otimes \bb C^n$ be the isometry given by
$Ve_j = e_j\otimes e_j$, where $\{e_j\}_{j=1}^n$ is the canonical basis of $\bb{C}^n$,
we have that
$\tilde{\Gamma}((x_{i,j})) = (V^* \otimes id_{\cl T})
\Gamma^{(n)}((x_{i,j}))(V \otimes id_{\cl T}).$
It now follows that
$\widetilde{\Gamma}$ is completely positive.
\end{proof}

We can now prove the main result of this section.

\begin{proposition}\label{opsysopsp} Let $\cl S$ and $\cl T$ be operator systems and let
$\tau$ be an operator system structure on $\cl S
  \otimes \cl T.$ Then the operator space $\cl S \otimes_{\tau} \cl
  T$  is an operator space tensor product of the operator spaces $\cl S$
  and $\cl T$ in the sense of \cite{bp}; that is, the following two conditions hold:
\begin{itemize}
\item[(1)] For any $s \in M_n(\cl S)$ and any $t \in M_m(\cl T)$ we have
$\|s \otimes t\|_{M_{mn}(\cl S \otimes_{\tau} \cl T)} \le
  \|s\|_{M_n(\cl S)}  \|t\|_{M_m(\cl T)}$.
\item[(2)] If $\phi: \cl S \to M_n$ and $\psi:\cl T \to M_m$ are completely
  bounded maps, then $\phi \otimes \psi: \cl S \otimes_{\tau} \cl T
  \to M_{mn}$ is completely bounded and $\|\phi \otimes \psi\|_{cb}
  \le \|\phi\|_{cb} \|\psi\|_{cb}$.
\end{itemize}
\end{proposition}
\begin{proof} Let $e$ denote the order unit of $\cl S$, and let $f$
  denote the order unit of $\cl T.$ To prove the first statement, it will be enough to assume that
  $\|s\| \le 1$ and $\|t\| \le 1$, and show that $\|s \otimes t\| \le
  1.$
But, in this case, $P=\begin{pmatrix} e_n & s\\s^* & e_n \end{pmatrix}
\in M_2(M_n(\cl S))^+ = M_{2n}(\cl S)^+$ and $Q= \begin{pmatrix} f_m &
  t\\t^* & f_m \end{pmatrix} \in M_2(M_m(\cl T))^+ = M_{2m}(\cl T)^+.$
Since $\tau$ is an operator system structure, Property (T2) implies that $P
\otimes Q \in M_{4mn}(\cl S \otimes_{\tau} \cl T)^+.$ Writing this
matrix in block form as a $4 \times 4$ matrix of $n \times m$ blocks, we have that
\[ \begin{pmatrix} e_n \otimes f_m & e_n \otimes t & s \otimes f_m & s
  \otimes t\\e_n \otimes t^* & e_n \otimes f_m & s \otimes t^* & s
  \otimes f_m\\s^* \otimes f_m & s^* \otimes t & e_n \otimes f_m & e_n
  \otimes t \\ s^* \otimes t^* & s^* \otimes f_m & e_n \otimes t^* &
  e_n \otimes f_m \end{pmatrix} \in M_4(M_{mn}(\cl S \otimes_{\tau}
\cl T))^+. \]
Compressing this block matrix to the four corner entries preserves positivity,
and hence
\[ \begin{pmatrix} e_n \otimes f_m & s \otimes t\\s^* \otimes t^* &
  e_n \otimes f_m \end{pmatrix} \in M_2(M_{mn}(\cl S \otimes_{\tau}
\cl T)), \]
and condition (1) follows.

To prove the second property, it will be enough to consider the case
where $\|\phi\|_{cb} \le 1$ and $\|\psi\|_{cb} \le 1.$  But in this
case, by \cite[Theorem~8.3]{Pa}, there exists a completely positive map $\Phi:M_2(\cl S) \to
M_2(M_n)$ given by
\[\Phi {\Huge (} \begin{pmatrix} s_{1,1} & s_{1,2}\\ s_{2,1} & s_{2,2} \end{pmatrix}
{\Huge )} = \begin{pmatrix} \phi_{1,1}(s_{1,1}) & \phi(s_{1,2}) \\ \phi(s_{2,1}^*)^* &
 \phi_{2,2}(s_{2,2}) \end{pmatrix} \in M_2(M_n) \] where
$\phi_{1,1},\phi_{2,2}: \cl S \to M_n$ are unital and completely
positive.  Also, there exists a similar completely positive map
$\Psi: M_2(\cl T) \to M_2(M_m)$ with analogous properties.

Let $\Phi_0 = \Phi \circ \delta: \cl S \to M_2(M_n)$ so that $\Phi_0(s) = \begin{pmatrix}
  \phi_{1,1}(s) & \phi(s) \\ \phi(s^*)^* &
  \phi_{2,2}(s) \end{pmatrix}$ and $\Psi_0: \cl T \to M_2(M_m)$ be defined in a
similar way. By Proposition \ref{tildetrick}, $\Phi_0$ and $\Psi_0$ are completely positive.
By Proposition \ref{non-unitalcp}, $\Phi_0 \otimes \Psi_0: \cl S \otimes_{\tau}
\cl T \to M_4(M_{mn})$ is completely positive.  Again, compressing to
corners yields a completely positive map $\Gamma: \cl S \otimes_{\tau}
\cl T \to M_2(M_{mn})$ with
\[\Gamma(s \otimes t) = \begin{pmatrix} \phi_{1,1}(s) \otimes
  \psi_{1,1}(t) & \phi(s) \otimes \psi(t) \\ \phi(s^*)^* \otimes
  \psi(t^*)^* & \phi_{2,2}(s) \otimes \psi_{2,2}(t) \end{pmatrix}.\]
Since $\phi\otimes\psi$ is a compression of a unital completely positive map,
it is completely contractive.
This completes the proof.
\end{proof}

One method that we shall use to distinguish operator system tensor products is to examine a canonical tensor product that they induce on the category of operator spaces and completely contractive maps.
Given an operator space $X$, there is a canonical operator
system $\cl S_X$ that can be associated to $X$.  If $X
\subseteq \cl B(H),$ then $\cl S_X \subseteq \cl B(H \oplus H)$ is the operator
system given by
$$\cl S_X = \left\{\left(
\begin{array}[c]{cc}
\lambda I_H & x\\
y^* & \mu I_H
\end{array}
\right) : \lambda,\mu\in \bb{C}, x, y\in X\right\}.$$ We regard $X
\subseteq \cl S_X,$ via the inclusion $x \to \begin{pmatrix} 0 & x\\0
  & 0 \end{pmatrix}.$ Note that the unit for $\cl S_X$ is $\begin{pmatrix} I_H & 0\\0 &I_H \end{pmatrix}.$

\begin{definition} \label{ind-op-space-def}
Let $X$ and $Y$ be operator spaces and $\tau$ be an operator system structure
on $\cl S_X \otimes \cl S_Y$. Then the embedding
$$X\otimes Y \subseteq \cl S_X\otimes_{\tau}\cl S_Y$$ endows $X\otimes Y$ with
an operator space structure; we call the resulting operator space
the \emph{induced operator space tensor product
of $X$ and $Y$} and denote it by $X\otimes^{\tau} Y$.
\end{definition}

\begin{proposition}\label{inducedopsp} Let $X$ and $Y$ be operator spaces,
let $\tau$ be an operator system structure on $\cl S_X \otimes \cl S_Y$,
and let $X \otimes^{\tau} Y$ be the induced operator space tensor product. Then $X \otimes^{\tau} Y$ is an operator space tensor product in the sense of \cite{bp}; that is, the following two conditions hold:
\begin{itemize}
\item[(1)] If $x \in M_n(X)$ and $y \in M_m(Y)$, then
$$\|x \otimes y\|_{M_{nm}(X \otimes^{\tau}Y)} \le \|x\|_{M_n(X)}\|y\|_{M_m(Y)}.$$
\item[(2)]  If $\phi:X \to M_n$ and $\psi:Y \to M_m$ are completely bounded, then
$\phi \otimes \psi: X \otimes^{\tau} Y \to M_{mn}$ is completely bounded and $\|\phi \otimes \psi\|_{cb}
  \le \|\phi\|_{cb} \|\psi\|_{cb}$.
\end{itemize}
\end{proposition}
\begin{proof}  The first claim follows from Proposition~\ref{opsysopsp} and the fact that the inclusions $X \subseteq \cl S_X$ and $Y \subseteq \cl S_Y$ are complete isometries.

To prove the second condition, note that by
\cite[Lemma 8.1]{Pa} if $\phi: X \to M_n$ is completely contractive, then the map $\Phi: \cl S_X \to M_2(M_n)$ given by
$$
\Phi \big( \begin{pmatrix} \lambda 1 & x_1\\ x_2^* & \mu 1 \end{pmatrix} \big) = \begin{pmatrix} \lambda I_n & \phi(x_1) \\ \phi(x_2)^* & \mu I_n \end{pmatrix}
$$
is a unital completely positive map. Similarly, the completely
contractive map $\psi:Y \to M_m$ yields a unital completely positive
map $\Psi: \cl S_Y \to M_2(M_m).$  By Property (T3) the map $\Phi
\otimes \Psi: \cl S_X \otimes_{\tau} \cl S_Y \to M_{4mn}$ is unital
and completely positive. Noticing that $\phi \otimes \psi$ occurs in a corner block of $\Phi \otimes \Psi$, we obtain that $\phi \otimes \psi$ is completely contractive.
\end{proof}

Let $\cl{OS}p$ be the the category whose objects are operator spaces and whose morphisms are
completely contractive linear maps.
Suppose that we are given an operator system tensor product $\tau: \cl O \times \cl O \to \cl O.$
We have that the mapping $\widetilde{\tau}
: \cl{OS}p\times\cl{OS}p\rightarrow \cl{OS}p$ given by $\widetilde{\tau}(X,Y) = X\otimes^{\tau} Y$
is an {\it operator space tensor product} in the sense of \cite{bp}.
We call $\widetilde{\tau}$ the operator space tensor product {\bf induced by} $\tau$.

The proof of the following result is similar to the proof of our last proposition, and we omit it.

\begin{proposition}\label{p_inhf}
If $\tau$ is a functorial operator system tensor product then $\widetilde{\tau}$ is a functorial operator space tensor product; that is, given any four operator spaces $X_1,X_2,Y_1$, and $Y_2$ and completely contractive maps $\phi:X_1 \to X_2$ and $\psi:Y_1 \to Y_2$, the map $\phi \otimes \psi: X_1 \otimes^{\tau} Y_1 \to X_2 \otimes^{\tau} Y_2,$ is completely contractive.
\end{proposition}


\section{The minimal tensor product}\label{s_min}

In this section we construct the operator system tensor product
$\min$, which is minimal among all operator system tensor products.
This section has overlaps with the work of Choi,
Effros and
Lance \cite{La1}, \cite{La}, \cite{Ef}, \cite{CEann}, \cite{CE2},
\cite{el} for C*-algebras and
Blecher and Paulsen \cite{bp} for operator spaces. We include this
material for completeness and because we will need some of the
results in later sections.

Let $\cl S$ and $\cl T$ be operator systems. For each $n\in \bb{N}$,
we let
\begin{align*} C_n^{\min} = C_n^{\min}(\cl S,\cl T) = \{(p_{i,j}) &\in M_n(\cl
S\otimes \cl T) : ((\phi\otimes\psi)(p_{i,j}))_{i,j} \in
M_{nkm}^+, \\
& \mbox{ for all } \phi\in S_k(\cl S), \psi\in S_m(\cl
T) \mbox{ for all } k,m\in \bb{N}\}.
\end{align*}

\begin{lemma} \label{l_gotoa}
Let $\cl S$ be an operator system and $P\in M_n(\cl S)$. If
$\phi^{(n)}(P)\in M_{nk}^+$ for every $\phi\in S_k(\cl S)$ and every
$k\in \bb{N}$, then $P\in M_n(\cl S)^+$.
\end{lemma}
\proof We may assume that $I\in \cl S\subseteq \cl B(H)$ for some Hilbert
space $H$. Suppose that $P = (p_{i,j})\in M_n(\cl S)$ and that
$\phi^{(n)}(P)\in M_{nk}^+$ for every $\phi\in S_k(\cl S)$ and every
$k\in \bb{N}$. Let $\xi = (\xi_1,\dots,\xi _n)^{\ttt}\in H^{n}$
(where $\ttt$ denotes transposition)
and $\phi : \cl S\rightarrow M_n$ be the mapping given by $\phi(x) =
((x\xi_j,\xi_i))_{i,j}$. We note that $\phi$ is completely positive.
Indeed, let $(x_{s,t})\in M_l(\cl S)^+$. We need to show that the
matrix $Y = (Y_{s,t})_{s,t}\in M_l(M_n)$, where $Y_{s,t} =
(x_{s,t}\xi_j,\xi_i)_{i,j}\in M_n$, is positive. Let $\lambda_s\in
\bb{C}^n$ for $s = 1,\dots,l$, where $\lambda_s =
(\lambda_{s,1},\dots,\lambda_{s,n})^{\ttt}$. Letting $\tilde{\lambda}
= (\lambda_1,\dots,\lambda_l)^{\ttt}$ and $\tilde{\xi}_s =
\sum_{i=1}^n \lambda_{s,i}\xi_i$, we have
\begin{eqnarray*}
(Y\tilde{\lambda},\tilde{\lambda}) & = & \sum_{s,t = 1}^l
(Y_{s,t}\lambda_t,\lambda_s) = \sum_{s,t = 1}^l \sum_{i,j=1}^n
(x_{s,t}\xi_j,\xi_i)\lambda_{t,j}\overline{\lambda_{s,i}}\\
& = & \sum_{s,t = 1}^l \left(x_{s,t}\left(\sum_{j=1}^n
\lambda_{t,j}\xi_j\right), \left(\sum_{i=1}^n
\lambda_{s,i}\xi_i\right)\right) = \sum_{s,t = 1}^l
\left(x_{s,t}\tilde{\xi}_t, \tilde{\xi}_s\right)\geq 0.
\end{eqnarray*}
Thus $\phi$ is completely positive and hence $\phi^{(n)}(P) =
(\phi(p_{i,j}))_{i,j}\in M_{n^2}^+$. Let
$$\eta = (e_1,e_2,\dots,e_n)^{\ttt}\in \bb{C}^{n^2},$$ where
$\{e_r\}_{r=1}^n$ is the standard basis of $\bb{C}^n$. Then
$$\sum_{i,j=1}^n (p_{i,j}\xi_j,\xi_i) = (\phi^{(n)}(P)\eta,\eta)\geq
0$$ and hence $P\in M_n(\cl B(H))^+$.
\endproof

In what follows we will identify $M_n(\cl S\otimes\cl T)$ with
$M_n(\cl S)\otimes\cl T$ in the natural way.

\begin{lemma}\label{l_slic}
Let $\cl S$ and $\cl T$ be operator systems and $P\in M_n(\cl
S)\otimes\cl T$. If $(\phi^{(n)}\otimes\psi)(P)\geq 0$ for all
$\phi\in S_{\infty}(\cl S)$ and all $\psi\in S_{\infty}(\cl T),$ then
$(\Phi\otimes\psi)(P)\geq 0$ for all $\Phi \in S_{\infty}(M_n(\cl
S))$ and all $\psi \in S_{\infty}(\cl T)$.
\end{lemma}
\proof Fix $m\in \bb{N}$ and $\psi\in S_m(\cl T).$ For each functional $\omega : M_m \to \bb C$,
let $\rho_{\omega} : M_n(\cl S)\otimes\cl T\rightarrow M_n(\cl S)$
be the mapping given by $\rho_{\omega}(X\otimes y) = \omega(\psi(y))X$, and
$L_{\omega} : M_m(V) \to V$ be the slice with respect to $\omega$. If $\eta_1,\eta_2\in \bb{C}^m$,
let $\omega_{\eta_1,\eta_2}$ be the functional on $M_m$ given by
$\omega_{\eta_1,\eta_2}(x) = (x\eta_1,\eta_2)$.

Suppose
that $(\phi^{(n)}\otimes\psi)(P) \in M_{nkm}^+$ for all $\phi\in
S_k(\cl S)$, $k \in \bb N$, and let
$\eta_1,\dots,\eta_r\in \bb{C}^m$. Since the map $(L_{\omega_{\eta_t,\eta_s}})_{s,t} : M_{nkm}
\to M_{nkr}$ is completely positive, we have that
$(L_{\omega_{\eta_t,\eta_s}}((\phi^{(n)}\otimes\psi)(P)))_{s,t}
\in M_{nkr}^+$.
Thus,
$$\phi^{(nr)}((\rho_{\omega_{\eta_t,\eta_s}}(P))_{s,t}) =
(\phi^{(n)}(\rho_{\omega_{\eta_t,\eta_s}}(P)))_{s,t}\geq 0, \ \ \mbox{ for all } \phi\in
S_{k}(\cl S), k\in \bb{N}.$$
By Lemma~\ref{l_gotoa}, $(\rho_{\omega_{\eta_t,\eta_s}}(P))_{s,t}\in
M_{nr}(\cl S)^+$, and hence $\Phi^{(r)}((\rho_{\omega_{\eta_t,\eta_s}}(P))_{s,t}) \geq 0$ for every
completely positive map $\Phi : M_n(\cl S) \rightarrow M_k$, $k\in \bb{N}$.
Fixing such a $\Phi$, we have that
$(L_{\omega_{\eta_t,\eta_s}}((\Phi\otimes\psi)(P)))_{s,t}\geq 0$.
Thus if $\xi_1,\dots,\xi_r\in \bb{C}^k$, then
$$\left((\Phi\otimes\psi)(P)\left(\sum_{t=1}^r \xi_t\otimes\eta_t\right),\left(\sum_{s=1}^r \xi_s\otimes\eta_s\right)\right) = $$
$$\left((L_{\omega_{\eta_t,\eta_s}}((\Phi\otimes\psi)(P)))_{s,t}(\xi_1,\dots,\xi_r)^{\ttt},(\xi_1,\dots,\xi_r)^{\ttt}\right) \geq 0.$$
It follows that $(\Phi\otimes\psi)(P)\geq 0$. The proof is complete.
\endproof

\begin{lemma}\label{l_phin}
If $\phi\in S_k(\cl S)$ and $\psi\in S_m(\cl T)$ then
$(\phi\otimes\psi)^{(n)} = \phi^{(n)}\otimes\psi$.
\end{lemma}
\proof It suffices to check the equality on elementary
tensors of the form $P = X\otimes y$, where $X = (x_{i,j})\in
M_n(\cl S)$ and $y\in \cl T$. For such a $P$ we have that
$(\phi^{(n)}\otimes\psi)(P) = (\phi(x_{i,j}))_{i,j}\otimes \psi(y)$.
On the other hand,
$$(\phi\otimes\psi)^{(n)}(P) = ((\phi\otimes\psi)(x_{i,j}\otimes y))_{i,j}
= (\phi(x_{i,j})\otimes\psi(y))_{i,j}.$$
\endproof

\begin{theorem}\label{th_spacial}
Let  $\cl S$ and $\cl T$ be operator systems, and let $\iota_{\cl S} :
\cl S\rightarrow \cl B(H)$ and $\iota_{\cl T} : \cl T\rightarrow \cl
B(K)$ be embeddings that are unital complete order isomorphisms onto their ranges.  The family
$(C_n^{\min}(\cl S,\cl T)_{n=1}^{\infty})$ is the operator system
structure on $\cl S\otimes\cl T$ arising from the embedding
$\iota_{\cl S}\otimes\iota_{\cl T} : \cl S\otimes\cl T\rightarrow
\cl B(H\otimes K)$.
\end{theorem}
\proof Let $P\in C_n^{\min}(\cl S,\cl T)$. We claim that
\begin{equation}\label{eq_ininc}
Q \stackrel{def}{=} (\iota_{\cl S}\otimes\iota_{\cl T})^{(n)}(P)\in \cl B((H\otimes K)^n)^+.
\end{equation}

Suppose that $Q = \sum_{r=1}^l X_r\otimes y_r$, where $X_r\in
M_n(\iota_{\cl S}(\cl S))$ and $y_r\in \iota_{\cl T}(\cl T)$ for $r = 1,\dots,l$. Let $\xi_s\in
H^{(n)}$ and $\eta_s\in K$ for $s = 1,\dots,k$, and set $\zeta =
\sum_{s=1}^k \xi_s\otimes\eta_s$. Let $\Phi : M_n(\iota_{\cl S}(\cl S))\rightarrow
M_k$ be the mapping given by $\Phi(X) = ((X\xi_t,\xi_s))_{s,t}$ and
let $\psi : \iota_{\cl T}(\cl T)\rightarrow M_k$ be the mapping given by $\psi(y) =
((y\eta_t,\eta_s))_{s,t}$. By the proof of Lemma~\ref{l_gotoa},
$\Phi$ and $\psi$ are completely positive. Since $Q\in
C_n^{\min}(\iota_{\cl S}(\cl S),\iota_{\cl T}(\cl T))$, Lemma~\ref{l_phin} implies that
$(\phi_0^{(n)}\otimes\psi_0)(Q)\in M_{nk^2}^+$, for all
$\phi_0\in S_k(\iota_{\cl S}(\cl S))$ and all $\psi_0\in S_k(\iota_{\cl T}(\cl T))$. Lemma~\ref{l_slic}
implies that $(\Phi\otimes\psi)(Q)\in M_{nk^2}^+$. Let $e =
(e_1,\dots,e_k)^{\ttt}\in \bb{C}^{k^2}$, where $\{e_j\}_{j=1}^k$
is the standard basis of $\bb{C}^k$. We then have
\begin{eqnarray*}
(Q\zeta,\zeta) & = & \sum_{r=1}^l \sum_{s,t=1}^k
(X_r\xi_t,\xi_s)(y\eta_t,\eta_s)\\
& = & \sum_{r=1}^l ((\Phi(X_r)\otimes\psi(y_r))e,e) =
((\Phi\otimes\psi)(Q)e,e).
\end{eqnarray*}
It follows that $Q\in \cl B((H\otimes K)^{n})^+$ and
(\ref{eq_ininc}) is established. Thus, if $D_n$ is the cone in
$M_n(\cl S\otimes\cl T)$ arising from the inclusion of $\iota_{\cl
S}(\cl S)\otimes\iota_{\cl T}(\cl T)$ into $\cl B(H\otimes K)$, we
have that $C_n^{\min}(\cl S,\cl T)\subseteq D_n$.

We now show that $D_n\subseteq C_n^{\min}(\cl S,\cl T)$. Suppose
that $\phi \in S_m (\cl S)$ and $\psi \in S_k(\cl T).$ By
identifying $\cl S = \iota_{\cl S}(\cl S) \subseteq \cl B(H)$ and
applying Arveson's extension theorem, we obtain a unital completely
positive map $\tilde{\phi}: \cl B(H) \to M_m$ that agrees with
$\phi$ on $\cl S.$ Similarly, we obtain a unital completely positive
map $\tilde{\psi}: \cl B(K) \to M_k$ that extends $\psi.$ By
C*-algebra theory, the minimal C*-tensor product
$\otimes_{\text{C*min}}$ satisfies $\cl B(H) \otimes_{\text{C*min}}
\cl B(K) \subseteq \cl B(H \otimes K)$ and there exists a unital
completely positive map $\tilde{\phi} \otimes \tilde{\psi}: \cl B(H)
\otimes_{\text{C*min}} \cl B(K) \to M_{mk}.$ Applying Arveson's
extension theorem once again, we obtain a unital completely positive
map $\gamma: \cl B(H \otimes K) \to M_{mk}.$ Therefore, if
$P=(p_{i,j}) \in D_n \subseteq \cl B((H \otimes K)^n)^+,$ then $(\phi
\otimes \psi(p_{i,j})) = (\gamma(p_{i,j})) \in M_{nmk}^+.$  Hence,
$D_n = C_n^{\min}(\cl S, \cl T).$

It follows that $C_n^{\min}(\cl S,\cl T)$ is an operator system
structure on $\cl S\otimes\cl T$ with an Archimedean matrix unit
$1\otimes 1$, where $1$ denotes the units for both $\cl S$ and $\cl
T$.
\endproof

\begin{definition}\label{d_min}
We call the operator system $(\cl S\otimes\cl T,(C_n^{\min}(\cl
S,\cl T))_{n=1}^{\infty}, 1 \otimes 1)$ the \emph{minimal tensor product} of $\cl S$ and $\cl T$
and denote it by $\cl S\otimes_{\min}\cl T$.
\end{definition}

\begin{theorem}\label{th_min}
The mapping $\min : \cl O\times\cl O\rightarrow\cl O$ sending $(\cl
S,\cl T)$ to $\cl S\otimes_{\min}\cl T$ is an injective, associative, symmetric, functorial
operator system tensor product.

Moreover, if $\cl S$ and $\cl T$ are operator systems and $\tau$ is
an operator system structure on $\cl S \otimes \cl T,$ then $\tau$
is larger than $\min$.
\end{theorem}

\begin{proof}
By Theorem \ref{th_spacial}, the mapping $\min$ is an injective
functorial operator system tensor product. Suppose that $\cl S_j$ is
an operator system and that $\iota_j: \cl S_j \to B(H_j)$ is a
complete order embedding, $j=1,2,3$. By the associativity of the
Hilbert space tensor product, we may identify $(H_1 \otimes H_2)
\otimes H_3$ with $H_1 \otimes (H_2 \otimes H_3).$ This
identification yields a complete order isomorphism of $(\cl S_1
\otimes_{\min} \cl S_2) \otimes_{\min} \cl S_3$ with $\cl S_1
\otimes_{\min} (\cl S_2 \otimes_{\min} \cl S_3)$, and hence $\min$
is associative. We see similarly that $\min$ is symmetric.

By (T3), we have that if $\tau$ is any operator
system structure on $\cl S\otimes\cl T$, then $M_n(\cl S \otimes_{\tau} \cl T)^+ \subseteq
C_n^{\min}(\cl S, \cl T)$ and hence $\min$ is the
minimal among all operator system structures on $\cl S\otimes \cl T$.
\end{proof}

\begin{remark}
{\rm It was shown in \cite{bp} that the minimal operator space tensor product, the spatial operator space tensor product, and the injective operator space tensor product all coincide.  For operator spaces $X$ and $Y$, we will let $X \check{\otimes} Y$ denote this tensor product, and choose to refer to it as the minimal operator space tensor product.}
\end{remark}

The following corollaries are immediate.

\begin{corollary}\label{p_restosp}
Let $X$ and $Y$ be operator spaces. Then the induced tensor product $X\otimes^{\min} Y$ (see Definition~\ref{ind-op-space-def}) coincides with the minimal operator space tensor product $X \check{\otimes} Y$.
\end{corollary}

\begin{corollary} \label{op-sp-min-cor}
Let $\cl S$ and $\cl T$ be operator systems. Then the identity map is a complete isometry
between the operator spaces
$\cl S \otimes_{\min} \cl T $ and $\cl S \check{\otimes} \cl T $.
\end{corollary}

\begin{corollary}\label{c_minC}
Let $A$ and $B$ be C*-algebras. Then the minimal operator system tensor product $A\otimes_{\min} B$ is completely order isomorphic to the image of $A\otimes B$ inside the minimal C*-algebraic tensor product $A \otimes_{\text{{\rm C*min}}} B$.
\end{corollary}

We close this section with a result which relates the minimal tensor product
of operator systems with the minimal operator system structure on an AOU space
studied in \cite{ptt}.
We recall from \cite{ptt} that if $(V,V^+)$ is an AOU space, $\OMIN(V)$ denotes the
minimal operator system whose underlying ordered $*$-vector space is $(V,V^+)$.

\begin{proposition}\label{p_ominmin}
Let $V$ and $W$ be AOU spaces. Equip the tensor product $V\otimes W$ with
the cone
$$Q_{\min} = \{u\in V\otimes W : (f\otimes g)(u) \geq 0, \mbox{ for all }
f\in S(V), g\in S(W)\}.$$ Then $\OMIN(V)\otimes_{\min}\OMIN(W) = \OMIN(V\otimes W)$.
\end{proposition}
\proof By \cite[Theorem 3.2]{ptt}, $\OMIN(V)\subseteq C(X)$, where
$X$ is the state space $S(V)$ equipped with the weak* topology.
Similarly, $\OMIN(W)\subseteq C(Y)$ where $Y = S(W)$. By the
injectivity of $\min$, we have that $\OMIN(V)\otimes_{\min}\OMIN(W)$
is an operator subsystem of $C(X)\otimes_{\min} C(Y)$. Denote the
 matrix
ordering on $\OMIN(V\otimes W)$ (respectively,
$\OMIN(V)\otimes_{\min}\OMIN(W)$) by $\{ Q_n \}_{n=1}^{\infty}$
(respectively, $\{ D_n\}_{n=1}^{\infty}$). Since $\OMIN(V\otimes W)$
is the minimal operator system structure on $(V\otimes W,Q_{\min})$,
we have that $D_n\subseteq Q_n$ for all $n\in \bb{N}$. Suppose that
$X = (x_{i,j})\in Q_n$. By \cite[Definition 3.1]{ptt},
$\sum_{i,j=1}^n \overline{\lambda_i}\lambda_jx_{i,j}\in Q_{\min}$
for all $\lambda_1,\dots,\lambda_n\in \bb{C}$. Thus, letting
$\tilde{\lambda} = (\lambda_1,\dots,\lambda_n)^{\ttt}$, we see that
for all $f\in S(V)$ and all $g\in S(W)$, we have
$$(((f\otimes g)(x_{i,j}))_{i,j}\tilde{\lambda},\tilde{\lambda}) = \sum_{i,j=1}^n
\overline{\lambda_i}\lambda_j (f\otimes g)(x_{i,j}) \geq 0.$$ It
follows that $X$ is a positive element of $M_n(C(X)\otimes_{\min}
C(Y)) \subseteq M_n(C(X\times Y))$, and hence $X\in D_n$.
Thus $D_n = Q_n$, for each $n\in \bb{N}$.
\endproof

\begin{remark}
{\rm Given two AOU spaces $V$ and $W$, which are also often called {\bf function systems,}
Effros \cite{Ef} (see also Namioka and Phelps \cite{np}) defines their minimal tensor product $V\otimes_{MIN} W$. The cone
$Q_{\min}$ from Proposition \ref{p_ominmin} coincides with the set of positive elements
of $V \otimes_{MIN} W$.  Thus, Proposition \ref{p_ominmin} says that
$\OMIN(V) \otimes_{\min} \OMIN(W) = \OMIN(V \otimes_{MIN} W).$}
\end{remark}


\section{The maximal tensor product}\label{s_max}

In this section we construct the maximal operator system tensor product and explore its properties.
Let $\cl S$ and $\cl T$ be operator systems whose units will both be denoted by 1. For each $n\in \bb{N}$,
we let
\begin{multline*} D_n^{\max} = D_n^{\max}(\cl S,\cl T) = \\ \{\alpha (P\otimes Q)\alpha^* :
P\in M_k(\cl S)^+, Q\in M_m(\cl T)^+, \alpha\in M_{n,km}, \ k,m\in \bb{N}\}.
\end{multline*}

\begin{lemma}\label{l_dmaxd}
Let $\cl S$ and $\cl T$ be operator systems and $\{D_n\}_{n=1}^{\infty}$ be
a compatible collection of cones, where $D_n\subseteq M_n(\cl S\otimes\cl T)$,
satisfying Property~(T2). Then
$D_n^{\max}\subseteq D_n$ for each $n\in \bb{N}$.
\end{lemma}
\begin{proof}
If $P\in M_k(\cl S)^+$ and $Q\in M_m(\cl T)^+$, Property~(T2)
implies that $P\otimes Q\in D_{km}$. The compatibility of
$\{ D_n\}_{n=1}^{\infty}$ implies that
$\alpha (P\otimes Q)\alpha^*\in D_n$ for every $\alpha\in M_{n,km}$.
Thus $D_n^{\max}\subseteq D_n$.
\end{proof}

\begin{lemma}\label{l_odot}
Let $\cl S$ and $\cl T$ be operator systems, $P = (P_{i,j})_{i,j}\in M_k(M_n(\cl S))^+$,
and $Q = (q_{i,j})_{i,j}\in M_k(\cl T)^+$. Then
$\sum_{i,j=1}^k P_{i,j}\otimes q_{i,j}\in D_n^{\max}.$
\end{lemma}
\proof
Let $I_n$ be the identity matrix in $M_n$, and
$X = (X_1,X_2,\dots,X_{k^2})\in M_{n,nk^2}$,
where $X_l\in M_n$ for $l = 1,\dots,k^2$, with
$$X_{1} = X_{k+2} = X_{2k+3} = \dots = X_{k^2} = I_n$$ and
$X_l = 0$ if $l\not\in \{1,k+2,2k+3,\dots,k^2\}$.
Then
$$\sum_{i,j=1}^k P_{i,j}\otimes q_{i,j} = X (P\otimes Q)X^*\in D_n^{\max}.$$
\endproof

\begin{proposition}\label{p_min}
Let $\cl S$ and $\cl T$ be operator systems. The family $\{
D_n^{\max}(\cl S,$ $\cl T)\}_{n=1}^{\infty}$ is a matrix ordering on
$\cl S\otimes\cl T$ with order unit $1\otimes 1$.
\end{proposition}
\begin{proof}
Let $n\in \bb{N}$.
Suppose that $\alpha_1 (P_1\otimes Q_1)\alpha_1^*$ and $\alpha_2 (P_2\otimes Q_2)\alpha_2^*$
are elements of $D_n^{\max}$, where
$P_i\in M_{k_i}(\cl S)^+, Q_i\in M_{m_i}(\cl T)^+$, and  $\alpha_i\in M_{n,k_i m_i}$ for $i = 1,2$.
Then
$\alpha_1 (P_1\otimes Q_1)\alpha_1^* + \alpha_2 (P_2\otimes Q_2)\alpha_2^*$ is equal to
$$
(\alpha_1,0,0,\alpha_2)((P_1\oplus P_2)\otimes (Q_1\oplus Q_2)) (\alpha_1,0,0,\alpha_2)^*,
$$
where $(\alpha_1,0,0,\alpha_2)\in M_{n,k_1m_1 + k_1m_2 + k_2m_1 + k_2m_2}$ and
$(P_1\oplus P_2)\otimes (Q_1\oplus Q_2)$ is identified with
$$
(P_1\otimes Q_1) \oplus (P_1\otimes Q_2) \oplus (P_2\otimes Q_1) \oplus (P_2\otimes Q_2).
$$
It is obvious that $D_n^{\max}$ is closed under positive scalar multiplies and
that $\{D_n^{\max}\}_{n=1}^{\infty}$ is a compatible family of cones.
By Lemma~\ref{l_dmaxd}, $D_n^{\max}\subseteq C_n^{\min}$, and hence
$D_n^{\max}\cap (-D_n^{\max}) \subseteq C_n^{\min} \cap (-C_n^{\min}) = \{0\}$.
Thus, $\{ D_n^{\max} \}_{n=1}^{\infty}$ is a matrix ordering. The fact that $1\otimes 1$ is
an order unit for $\{D_n^{\max}\}_{n=1}^{\infty}$ follows from the inclusions
$D_n^{\max}\subseteq C_n^{\min}$ and the fact that it is a matrix order unit for
$\{C_n^{\min}\}_{n=1}^{\infty}$.
\end{proof}

\begin{definition}\label{d_max}
Let $C_n^{\max} = C_n^{\max}(\cl S,\cl T)$ be the
Archimedeanization of the matrix ordering $\{ D_n^{\max}(\cl S,\cl
T) \}_{n=1}^{\infty}$. We call the operator system $$(\cl S\otimes\cl
T, \{C_n^{\max}(\cl S,\cl T)\}_{n=1}^{\infty}, 1 \otimes 1)$$ the \emph{maximal operator system
tensor product} of $\cl S$ and $\cl T$ and denote it by $\cl
S\otimes_{\max}\cl T$.
\end{definition}

By \cite[Remark 3.19]{ptt}, we have that $P \in C_n^{\max}(\cl S,
\cl T)$ if and only if $re_n +P \in D_n^{\max}(\cl S, \cl T)$ for
every $r>0.$

\begin{theorem}\label{th_max}
The mapping $\max : \cl O\times\cl O\rightarrow\cl O$ sending $(\cl
S,\cl T)$ to $\cl S\otimes_{\max}\cl T$ is a symmetric, associative, functorial
operator system tensor product.  Moreover, if $\tau$ is an operator system structure on $\cl S \otimes \cl T,$ then $\max$ is larger than $\tau$.
\end{theorem}
\begin{proof}
Let $\cl S$ and $\cl T$ be operator systems. By its definition,
the family $\{C_n^{\max}\}_{n=1}^{\infty}$ satisfies Property~(T1) and Property~(T2). Since
$C_n^{\max}(\cl S,\cl T)\subseteq C_n^{\min}(\cl S,\cl T)$, it
follows from Theorem \ref{th_min} that $\cl S\otimes_{\max}\cl T$ satisfies Property~(T3).
Suppose that $\phi\in \operatorname{UCP}(\cl S_1,\cl S_2)$ and $\psi\in \operatorname{UCP}(\cl T_1,\cl T_2)$,
and let $P\in M_k(\cl S_1)^+$, $Q\in M_m(\cl T_1)^+$, and $\alpha\in M_{n,km}$.
Then $\phi^{(k)}(P)\in M_k(\cl S_2)^+$ and $\psi^{(m)}(Q)\in M_m(\cl T_2)^+$.
Hence
$$(\phi\otimes\psi)^{(n)}(\alpha (P\otimes Q)\alpha^*) =
\alpha (\phi^{(k)}(P)\otimes \psi^{(m)}(Q))\alpha^* \in M_n(\cl S_2\otimes_{\max}\cl T_2)^+.$$
It follows that $(\phi\otimes\psi)^{(n)}(D_n^{\max}(\cl S_1,\cl T_1))\subseteq D_n^{\max}(\cl S_2,\cl T_2)$.
Lemma~\ref{l-cparch} now implies that Property~(T4) is satisfied.

Suppose that $P\in M_k(\cl S)^+$ and $Q\in M_m(\cl T)^+$.
Recall that the map $\theta: \cl S\otimes\cl T \rightarrow\cl T\otimes\cl S$ is given by
$\theta(x\otimes y) = y\otimes x$. We have that, after conjugation with a permutation matrix,
$\theta^{(km)}(P\otimes Q) = Q\otimes P$. It follows that if $\alpha\in M_{n,km}$, then
$$\theta^{(n)}(\alpha(P\otimes Q)\alpha^*) = \alpha \theta^{(km)}(P\otimes Q)\alpha^* =
\alpha (Q\otimes P)\alpha^*.$$
Thus $\theta : \cl S\otimes_{\max}\cl T\rightarrow \cl T\otimes_{\max}\cl S$ is a complete order isomorphism
and hence $\max$ is symmetric.

The fact that max is the maximal operator system
tensor product follows from Lemma~\ref{l_dmaxd}. It remains to prove associativity.
Let $\cl R, \cl S,$ and $\cl T$ be operator systems. The inclusion $p\rightarrow p\otimes 1$ of
$\cl R \otimes \cl S$ into
$\cl R \otimes_{\max} ( \cl S \otimes_{\max} \cl T)$
endows $\cl R \otimes \cl S$ with an operator system structure and, by the maximality
of $\max$, it
yields a completely positive map $\gamma: \cl R \otimes_{\max}
\cl S \to \cl R \otimes_{\max} ( \cl S \otimes_{\max} \cl T).$ If
$s: \cl T \to \bb C$ is any state, then by functoriality there exists a
completely positive map $\id_{\cl S} \otimes_{\max} s: \cl S
\otimes_{\max} \cl T \to \cl S \otimes_{\max} \bb C = \cl S.$
Functoriality also gives a completely positive map $\id_{\cl R}
\otimes_{\max} (\id_{\cl S} \otimes_{\max} s): \cl R \otimes_{\max} (
\cl S \otimes_{\max} \cl T) \to \cl R \otimes_{\max} \cl S$ that is
easily seen to be a left inverse for $\gamma.$  Hence $\gamma$ is a
complete order isomorphism onto its range.
Let $\gamma_1 : \gamma(\cl R\otimes_{\max}\cl S)\times \cl T \rightarrow \cl R\otimes_{\max}
(\cl S\otimes_{\max}\cl T)$ be the map sending $(p\otimes 1, z)$ to $p\otimes z$
and $\tilde{\gamma}_1 :
(\cl R \otimes_{\max} \cl S) \otimes \cl T \to \cl R \otimes_{\max} (\cl S \otimes_{\max} \cl T)$
be the corresponding linear map.
The map $\tilde{\gamma}_1$ endows
$(\cl R \otimes_{\max} \cl S) \otimes \cl T$ with an operator system structure.
It follows that $\tilde{\gamma}_1$ is
completely positive from
$(\cl R \otimes_{\max} \cl S) \otimes_{\max} \cl T$ to
$\cl R \otimes_{\max} ( \cl S \otimes_{\max} \cl T).$ However,
$\tilde{\gamma}_1$ coincides with the canonical mapping
from $(\cl R \otimes\cl S)\otimes \cl T$ onto $\cl R \otimes (\cl S\otimes \cl T)$.
Thus, the matricial cones of $(\cl R \otimes_{\max} \cl S) \otimes_{\max} \cl T$ are
contained in the corresponding matricial cones of $\cl R \otimes_{\max} ( \cl S \otimes_{\max} \cl T).$
A similar argument shows the converse inclusions, and hence we have that
$(\cl R \otimes_{\max} \cl S) \otimes_{\max} \cl T = \cl R \otimes_{\max} ( \cl S \otimes_{\max} \cl T).$
\end{proof}

\begin{definition}\label{d_jcp}
Let $\cl S$ and $\cl T$ be operator systems. A bilinear map $\phi :
\cl S\times\cl T\rightarrow \cl B(H)$ is called
\emph{jointly completely positive} if $\phi^{(n,m)}(P,Q)$ is a positive element of $M_{nm}(\cl B(H))$,
for all $P\in M_n(\cl S)^+$ and all $Q\in M_m(\cl T)^+$.
\end{definition}

The following result from \cite{La1} gives a useful characterization
 of jointly completely
positive maps. Given a bounded bilinear map $\phi: \cl S \times \cl
T \to \bb C$ we can define $\cl L(\phi): \cl S \to \cl T^d$
(respectively, $\cl R(\phi): \cl T \to \cl S^d$) by $\cl
L(\phi)(s)(t) = \phi(s,t)$ (respectively, $\cl R(\phi)(t)(s) =
\phi(s,t)$).

\begin{lemma}\label{llance} \textnormal{(\cite[Lemma~3.2]{La1})\textbf{.}}
Let $\cl S$ and $\cl T$ be operator systems and
let
  $\phi: \cl S \times \cl T \to \bb C$ be a bilinear map.  Then the following
  are equivalent:
\begin{itemize}
\item[(i)] $\phi$ is jointly completely positive.
\item[(ii)] $\cl L(\phi): \cl S \to \cl T^d$ is completely positive.
\item[(iii)] $\cl R(\phi): \cl T \to \cl S^d$ is completely positive.
\end{itemize}
\end{lemma}


The next theorem characterizes the maximal operator system tensor product in terms of a certain
universal property.

\begin{theorem}\label{th_linjcp}
Let $\cl S$ and $\cl T$ be operator systems.

\begin{itemize}
\item[(i)] If $\phi : \cl S\times\cl T\rightarrow \cl B(H)$ is a jointly completely positive map,
then its linearization $\phi_L : \cl S\otimes\cl T\rightarrow \cl B(H)$
is completely positive on $\cl S\otimes_{\max}\cl T$.

\item[(ii)] If $\psi : \cl S\otimes_{\max}\cl T\rightarrow \cl B(H)$ is completely positive, then
the map $\phi : \cl S\times \cl T\rightarrow \cl B(H)$ given by $\phi(x,y) = \psi(x\otimes y)$,
for $x\in \cl S$ and $y\in \cl T$, is jointly completely positive.

\item[(iii)] If $\tau$ is an operator system structure on $\cl S\otimes\cl T$
with the property that the linearization of every unital jointly completely positive map $\phi : \cl S\times\cl T\rightarrow
\cl B(H)$ is completely positive on $\cl S\otimes_{\tau}\cl T$, then
$\cl S\otimes_{\tau} \cl T = \cl S\otimes_{\max}\cl T$.

\item[(iv)] For every $n \in \bb N,$ we have that
$$C_n^{\max}(\cl S, \cl T) =
\{u \in M_n(\cl S \otimes \cl T) : \phi_L^{(n)}(u) \ge 0, \ \ \text{
for all jointly completely }$$
$$\mbox{positive } \phi: \cl S \times \cl T \to \cl B(H) \text{ and all Hilbert spaces } H\}.$$
\end{itemize}
\end{theorem}
\begin{proof}
Fix operator systems $\cl S$ and $\cl T$.

(i) Let $\phi : \cl S\times\cl T\rightarrow \cl B(H)$ be a jointly completely positive map.
If $P\in M_k(\cl S)^+$ and $Q\in M_m(\cl T)^+$, then
$\phi_L^{(km)}(P\otimes Q) = \phi^{(k,m)}(P,Q)\geq 0$. Thus if $\alpha\in M_{n,km}$, then
$$\phi_L^{(n)}(\alpha (P\otimes Q) \alpha^*) = \alpha \phi_L^{(km)}(P\otimes Q)\alpha^* \geq 0,$$
and hence $\phi_L^{(n)}(D_n^{\max}) \subseteq M_n(\cl B(H))^+$. By Lemma~\ref{l-cparch}, we have
$\phi_L$ is completely positive.

(ii) If $P\in M_k(\cl S)^+$ and $Q\in M_m(\cl T)^+$, then
$\phi^{(k,m)}(P,Q) = \psi^{(km)}(P\otimes Q) \geq 0$.

(iii) By Lemma~\ref{l_dmaxd}, $\max$ is larger than $\tau$, and hence
every unital completely positive map on $\cl S\otimes_{\tau}\cl T$ is
completely positive on $\cl S\otimes_{\max}\cl T$. By hypothesis,
$\operatorname{UCP}(\cl S\otimes_{\tau}\cl T,\cl B(H)) = \operatorname{UCP}(\cl S\otimes_{\max}\cl T,\cl B(H))$
for every Hilbert space $H$. By Lemma~\ref{l-detoss}, we have $\cl S\otimes_{\tau}\cl T =
\cl S\otimes_{\max}\cl T$.

(iv) Let $C_n \subseteq M_n(\cl S \otimes \cl T)$ be the set defined
by the right hand side of the displayed equation, and check that $\{ C_n
\}_{n=1}^{\infty}$ is an operator system structure, say $\tau$, on $\cl S
\otimes \cl T$. The result now follows
by observing that $\tau$ satisfies the hypotheses of (iii).
\end{proof}

If $X$ and $Y$ are operator spaces, then we let
$X\hat{\otimes} Y$ denote the operator space projective tensor product.
We refer the reader to \cite{bp} and \cite{er} for the definition and properties of this tensor product.

\begin{theorem}\label{th_resmaxi}
Let $X$ and $Y$ be operator spaces. Then $X\otimes^{\max} Y$ coincides with the
operator space projective tensor product $X\hat{\otimes} Y$.
\end{theorem}

\begin{proof}  Let $e= \begin{pmatrix} e_1 & 0\\0 & e_2 \end{pmatrix}$ denote the identity of $\cl S_X$ and let $f= \begin{pmatrix} f_1 & 0\\0 & f_2 \end{pmatrix}$ denote the identity of $\cl S_Y$, so that $e \otimes f$ is the identity of $\cl S_X \otimes \cl S_Y.$ Let $U=(u_{r,s}) \in M_p(X\otimes^{\max} Y)$ with $\|U\|^{\max} < 1.$ We must prove that the norm
$\|U\|$ of $U$ as an element of $M_p(X \hat{\otimes} Y)$ does not exceed 1.

We have that
\begin{align*}
\begin{pmatrix} \|U\|^{\max}(e\otimes f)_p & U\\ U^* & \|U\|^{\max}(e \otimes f)_p \end{pmatrix} \in& \ M_2(M_p(\cl S_X \otimes_{\max} \cl S_Y))^+ \\
&= C_{2p}^{\max}(\cl S_X, \cl S_Y)
\end{align*}
and hence
\begin{multline*} \begin{pmatrix} (e \otimes f)_p & U\\ U^* & (e \otimes f)_p \end{pmatrix} = \\ (1- \|U\|^{\max}) \begin{pmatrix} (e\otimes f)_p & 0\\0 & (e\otimes f)_p \end{pmatrix} +
\begin{pmatrix} \|U\|^{\max}(e\otimes f)_p & U\\ U^* & \|U\|^{\max}(e \otimes f)_p \end{pmatrix}\end{multline*}
is in $D_{2p}^{\max}(\cl S_X, \cl S_Y).$

Thus, there exist $P=(P_{i,j}) \in M_n(\cl S_X)^+, Q=(Q_{i,j}) \in M_m(\cl S_Y)^+$ and a $2p \times mn$ matrix $T= \begin{pmatrix} A\\B \end{pmatrix}$ where $A= (a_{r,(i,k)}), B= (b_{r,(i,k)})$ are $p \times mn$ matrices, such that
\[ \begin{pmatrix} (e \otimes f)_p & U\\ U^* & (e \otimes f)_p \end{pmatrix} = T(P\otimes Q)T^*.\]
This leads to the equations
$(e \otimes f)_p = A(P \otimes Q)A^*,$  $U= A(P \otimes Q)B^*,$ $U^* = B(P\otimes Q)A^*,$ and $(e \otimes f)_p = B(P\otimes Q)B^*.$

Recall that each element of $\cl S_X$ and $\cl S_Y$ is itself a $2
\times 2$ matrix and let $P_{i,j}= \begin{pmatrix} \alpha_{i,j}e_1 &
x_{i,j}\\w_{i,j}^* & \beta_{i,j}e_2\end{pmatrix} \in \cl S_X,$ where
$\alpha_{i,j}, \beta_{i,j} \in \bb C$ and $x_{i,j}, w_{i,j} \in X.$
Similarly, let $Q_{k,l} = \begin{pmatrix} \gamma_{k,l}f_1 & y_{k,l}
\\ z_{k,l}^* & \delta_{k,l}f_2 \end{pmatrix} \in \cl S_Y,$ where
$\gamma_{k,l}, \delta_{k,l} \in \bb C$ and $y_{k,l}, z_{k,l} \in Y.$
Finally, set $R_1= ( \alpha_{i,j})$, $R_2=(\beta_{i,j})$,
$S_1=(\gamma_{k,l})$, $S_2= (\delta_{k,l})$, $\cl X= (x_{i,j})$, and
$\cl Y = (y_{k,l}).$

Since $P$ and $Q$ are positive we have that $R_1$, $R_2$, $S_1$, and $S_2$
are positive scalar matrices, that $(w_{i,j}^*) = \cl X^*$, $(z_{k,l}^*) = \cl Y^*$, and that
for every $r>0,$
$\|(R_1 + rI_n)^{-1/2} \cl X(R_2 + rI_n)^{-1/2}\| \le 1$ in $M_n(X)$
and $\|(S_1 +rI_m)^{-1/2}\cl Y(S_2 + rI_m)^{-1/2}\| \le 1$ in $M_m(Y)$
(see \cite[p. 99]{Pa}).

Let $R_1e_1$ denote the matrix $(\alpha_{i,j} e_1)$ with similar
definitions for $R_2e_2,$ $S_1f_1$, $S_2f_2.$ Recalling that the
equation $(e \otimes f)_p = A(P \otimes Q)A^*$ takes place in $\cl
S_X \otimes \cl S_Y,$ which is represented by $4 \times 4$ block
matrices, we see that it yields $(e_i \otimes f_j)_p = A(R_ie_i
\otimes S_jf_j)A^*$ for $i,j =1,2.$  Thus, $I_p = A(R_i \otimes
S_j)A^*.$ Similarly, $I_p= B(R_i \otimes S_j)B^*.$

Recall that we have identified $x$ with $\begin{pmatrix} 0 & x\\0 & 0 \end{pmatrix}$ and $y$ with $\begin{pmatrix} 0 & y\\ 0 & 0 \end{pmatrix}$, so that $U$ only occurs in the $(1,4)$ block of the $4 \times 4$ block matrix, with the remaining entries equal to zero.
Thus, the equation $U = A(P\otimes Q)B^*$ in $\cl S_X \otimes \cl S_Y$ yields $U=A(\cl X \otimes \cl Y)B^*$ in $X \otimes Y.$

In the case that all scalar matrices $R_1,R_2,S_1$ and $S_2$ are
invertible, let $A_1 = A(R_1 \otimes S_1)^{1/2}$ and let $B_1 =
B(R_2 \otimes S_2)^{1/2},$ so that $U = A_1(R_1 \otimes S_1)^{-1/2}(
\cl X \otimes \cl Y)(R_2 \otimes S_2)^{-1/2}B_1^* =
A_1[(R_1^{-1/2}\cl XR_2^{-1/2}) \otimes (S_1^{-1/2} \cl Y
S_2^{-1/2})] B_1^*.$ Since $A_1A_1^* = I_p$ and $B_1B_1^* = I_p$, we
have that $\|R_1^{-1/2} \cl X R_2^{-1/2}\| \le 1$ and $\|S_1^{-1/2}
\cl Y S_2^{-1/2}\| \le 1,$ and we have obtained $U = A_1( \cl X_1
\otimes \cl Y_1)B_1^*$, where
$\cl X_1 = R_1^{-1/2}\cl XR_2^{-1/2}$, $\cl Y_1 = S_1^{-1/2} \cl Y
S_2^{-1/2}$ and all matrices $A_1,\cl X_1, \cl Y_1, \cl B_1$ have
norm at most one. This implies that $\|U\| \le 1$.

When the scalar matrices are not all invertible, one needs to first
add $rI_n$ and $rI_m$ ($r > 0$) to the corresponding matrices, set
$A_1 = A[(R_1 + rI_n) \otimes (S_1 +rI_m)]^{1/2}$, $B_1 = B[(R_2 +
rI_n) \otimes (S_2 + rI_m)]^{1/2}$, and conclude that $\|U\| \le 1 +
Cr$ where $C$ is a constant independent of $r.$ Since this
inequality holds for all $r > 0$, we again obtain that $\|U\|\leq
1$.
\end{proof}

\begin{remark}\label{rem-fdopsys}
{\rm Given two operator systems $\cl S$ and $\cl T$, Choi and Effros
define in \cite{CEann} an ordered $*$-vector space, which they call
the maximal tensor product of $\cl S$ and $\cl T$, using a scalar
version of Theorem~\ref{th_linjcp} (iv) to define its positive cone.
Let $A$ and $B$ be C*-algebras. Then $C_n^{\min}(A,B)$ can be
canonically identified with $C_1^{\min}(M_n(A),B)$ and any bilinear
map $\phi: M_n(A) \times B \to \bb{C}$ can be identified with a
bilinear map $\tilde{\phi} : A \times B \to M_n$. Using techniques
of Lance \cite{La} and these identifications, one can show that $u
\in C_n^{\min}(A,B)$ if and only if $\phi_L^{(n)}(u) \ge 0$ for all
$H$ and for all $\phi: A \times B \to \cl B(H)$ with $\phi$ jointly
completely positive and of finite rank. (We say that a bounded
bilinear map $\phi: A \times B \to \cl B(H)$ is of finite rank if
the induced map $\cl L(\phi) : A \to \cl B(B,\cl B(H))$ has finite
rank.) This fails for general operator systems, as we shall now
show. If $\cl S$ is a finite-dimensional operator system, then for
any operator system $\cl T$, every bilinear map $\phi: \cl S \times
\cl T \to \cl B(H)$ is of finite rank. Thus, if Lance's result held
for operator systems, it would imply that the minimal and maximal
tensor products on $\cl S \otimes \cl T$ are equal whenever $\cl S$
is finite dimensional. Applying this fact to operator systems of the
form $\cl S_X$ and using Corollary~\ref{op-sp-min-cor} and
Theorem~\ref{th_resmaxi} would yield that $X \hat\otimes Y$ is
completely isometric to $X \otimes_{\min} Y$ whenever $X$ is a
finite-dimensional operator space.  But this is known to be false,
see \cite{bp}. Thus, the analogue of this result of Lance fails for
operator systems. In particular, we see that there exist
finite-dimensional operator systems that are not
$(\min,\max)$-nuclear. Thus, the characterization due to \cite{Ki} and \cite{CE3} of nuclearity of C*-algebras
via the completely positive approximatation property(CPAP) does not hold
for operator systems.}
\end{remark}

Even for matrix algebras, the maximal operator space cross-norm is
larger than the operator space norm induced by the maximal operator
system tensor product. In fact, it can be shown that the cb-norm of
$\id: M_n \otimes_{\max} M_n \to M_n \hat\otimes M_n$ tends to
$+\infty$ as $n \to +\infty.$ One way to prove this is to use
Theorem~\ref{th_maxcst} below to see that $M_n \otimes_{\max} M_n =
M_{n^2},$ up to a unital complete order isomorphism, use the fact that the norm
on $M_n \hat\otimes M_n$ is larger than the Haagerup tensor norm
\cite{bp} and compare these two norms for the element $U =
\sum_{i=1}^n E_{1,i} \otimes E_{i,1}.$

The following result characterizes when these two tensor products yield completely isomorphic operator spaces.

\begin{proposition}\label{p_maxeq}
Let $\cl S$ and $\cl T$ be operator systems. The following are equivalent:
\begin{itemize}
\item[(i)] \ The identity map $\psi: \cl S\otimes_{\max}\cl T \rightarrow \cl S\hat{\otimes}\cl T$ is completely bounded.
\item[(ii)] There exists $C > 0$ such that for every jointly completely contractive map
$\phi : \cl S\times\cl T\rightarrow \cl B(H)$ there exist jointly completely
positive maps $\phi_i : \cl S\times\cl T\rightarrow \cl B(H)$
such that $\|\phi_i(e_{\cl S},e_{\cl T})\|\leq C$, $i = 1,2,3,4$,
and $\phi = (\phi_1 - \phi_2) + i(\phi_3 - \phi_4)$.
\end{itemize}
\end{proposition}

\begin{proof}(i)$\Rightarrow$(ii). By assumption, the identity map $\psi : \cl S\otimes_{\max}\cl T\rightarrow\cl S\hat{\otimes}\cl T$
is completely bounded; let $C$ be its cb-norm. Let $\phi : \cl S\times\cl T\rightarrow \cl B(H)$ be
a jointly completely contractive map. Then its linearization
$\tilde{\phi} : \cl S\hat{\otimes}\cl T \rightarrow \cl B(H)$ is completely contractive and hence
$\tilde{\phi}\circ\psi : \cl S\otimes_{\max}\cl T \rightarrow\cl B(H)$ is
completely bounded with cb-norm not exceeding $C$. By the Wittstock Decomposition Theorem,
there exist completely positive maps $\tilde{\phi}_i : \cl S\otimes_{\max}\cl T\rightarrow\cl B(H)$ for
$i = 1,2,3,4$, with norm not exceeding $C$
and such that
$\tilde{\phi} = (\tilde{\phi}_1 - \tilde{\phi}_2) + i(\tilde{\phi}_3 - \tilde{\phi}_4)$.
If $\phi_i$ is the bilinear map corresponding to $\tilde{\phi}_i$
then $\phi_i$ ($i = 1,2,3,4$) is jointly completely positive by Theorem \ref{th_linjcp}(ii); clearly,
$\phi = (\phi_1 - \phi_2) + i(\phi_3 - \phi_4)$.

(ii)$\Rightarrow$(i). Let $\iota : \cl S\hat{\otimes}\cl T \rightarrow \cl B(H)$ be a
complete isometry. By assumption, $\iota = (\tilde{\phi}_1 - \tilde{\phi}_2) + i(\tilde{\phi}_3 - \tilde{\phi}_4)$,
where $\tilde{\phi}_i$ is the linearization of a jointly completely positive map $\phi_i : \cl S\times\cl T\rightarrow\cl B(H)$ for $i = 1,2,3,4$.
By Theorem \ref{th_linjcp}(i), $\tilde{\phi}_i : \cl S\otimes_{\max}\cl T\rightarrow\cl B(H)$
is completely positive, and hence completely bounded. It follows that the identity map
$\id: \cl S\otimes_{\max}\cl T\rightarrow \cl B(H)$ is completely bounded, and therefore
$\cl S\otimes_{\max}\cl T$ is completely boundedly isomorphic to $\cl S\hat{\otimes}\cl T$.
\end{proof}

Except for the last conclusion, the following result is a consequence of the deep work of Choi, Effros, and Lance (see \cite{CEann}, \cite{CE2}, \cite{CE3}, and \cite{el}).

\begin{theorem}\label{th_maxcst}
Let $A$ and $B$ be C*-algebras. Then the operator system $A\otimes_{\max} B$ is
completely order isomorphic to the image of $A\otimes B$ inside the maximal
C*-algebraic tensor product of $A$ and $B$.
\end{theorem}
\begin{proof}
Let $\cl C = A  \otimes_{\text{C*max}} B$ denote the maximal C*-algebraic tensor product of $A$ and $B$. We claim that the faithful inclusion $A\otimes B \subseteq \cl C$ endows $A \otimes B$ with an operator system structure. Indeed, (T1) and (T2) are trivial and (T3) follows since it holds for the minimal C*-tensor product, which is a quotient of $\cl C$. We let $A\otimes_{\tau}B \subseteq \cl C$ denote this operator system.

For each $n\in \bb{N}$, let $D_n = M_n(A\otimes_{\tau} B)^+ =
M_n(A\otimes B) \cap M_n(\cl C)^+$.
Lemma \ref{l_dmaxd} implies that $A\otimes_{\max} B$ is larger than
$A \otimes_{\tau} B,$ and hence $C_n^{\max}(A,B) \subseteq D_n.$

We next show that the AOU spaces $(M_n(A\otimes B), C_n^{\max}(A,B))$ and
$(M_n(A\otimes B), D_n)$ have the same state space. In view of the last inclusion,
it suffices to show that if $f : A\otimes B \rightarrow \bb{C}$ and $f(C_n^{\max}(A,B))\subseteq \bb{R}^+$ then
$f(D_n) \subseteq \bb{R}^+$. So, let us fix an $f$ with $f(C_n^{\max}(A,B))\subseteq \bb{R}^+$.
Suppose that
$X = \sum_{i=1}^k a_i\otimes b_i$,
with $a_i\in M_n(A)$ and $b_i\in B$. Then
$$XX^* = \sum_{i,j=1}^k a_ia_j^*\otimes b_ib_j^*.$$
Let $P = (a_ia_j^*)_{i,j}$ and $Q = (b_ib_j^*)_{i,j}$; then $P\in M_k(M_n(A))^+$ and $Q\in M_k(B)^+$.
It follows from Lemma \ref{l_odot} that
$XX^*\in C_n^{\max}(A,B)$ and hence $f(XX^*) \geq 0$.
On the other hand, by the associativity of the C*-algebraic tensor product and the fact that
$M_n$ is a nuclear C*-algebra, we have a natural identification
$M_n(\cl C) = M_n(A)\otimes_{\text{C*max}} B$.
By the definition of the set of states on the C*-algebraic tensor product \cite[p. 381]{La},
we have that $f(D_n)\subseteq \bb{R}^+$.

Now let $u\in D_n$ and $f : M_n(A\otimes_{\max} B)\rightarrow \bb{C}$ be
positive, that is, $f(C_n^{\max}(A,B)) \subseteq \bb{R}^+$.
By the previous paragraph, $f(u) \geq 0$. By \cite[Proposition 3.13]{pt}, $u\in C_n^{\max}(A,B)$ and
the proof is complete.
\end{proof}

For the next proposition, we recall that if $(V,V^+)$ is an AOU space, $\OMAX(V)$ denotes the
maximal operator system whose underlying ordered $*$-vector space is $(V,V^+)$ \cite{ptt}.

\begin{proposition}\label{p_omaxmax}
Let $(V,V^+)$ and $(W,W^+)$ be AOU spaces. Equip the tensor product
$V\otimes W$ with the Archimedenization of the cone
$$Q_{\max} = \left\{\sum_{i=1}^k v_i\otimes w_i : v_i\in V^+, w_i\in W^+,  \text{ and } k\in\bb{N} \right\}.$$
Then $\OMAX(V)\otimes_{\max}\OMAX(W) = \OMAX(V\otimes W)$.
\end{proposition}
\proof
Recall that the matrix ordering on $\OMAX(V)$ is the
Archimedeanization of $\{ D_n^{\max}(V) \}_{n=1}^{\infty}$ where
$$D_n^{\max}(V) = \left\{\sum_{j=1}^k a_j\otimes v_j : a_j\in M_n^+, v_j\in V^+,  \text{ and } k\in\bb{N} \right\}.$$
Define similarly $\{D_n^{\max}(W)\}_{n=1}^{\infty}$ with respect to the cone $W^+$
and $\{D_n^{\max}(V\otimes W)\}_{n=1}^{\infty}$ with respect to the cone $Q_{\max}$.
It suffices to show that
$$D_n^{\max}(V\otimes W) = \left\{\alpha (P\otimes Q) \alpha^* : P\in D_k^{\max}(V), Q\in D_m^{\max}(W),
\alpha\in M_{n,km}\right\}.$$
Let $D_n$ denote the right hand side of the last equation.
If $a_j\in M_n^+$ and $\sum_{i=1}^{k_j} v_i^j\otimes w_i^j \in Q_{\max}$, $j = 1,\dots,l$, where $v_i^j\in V^+$ and
$w_i^j\in W^+$, then
$$\sum_{j=1}^l a_j\otimes \left(\sum_{i=1}^{k_j} v_i^j\otimes w_i^j\right) =
\sum_{j,i} a_j\otimes v_i^j\otimes w_i^j.$$
Since $\sum_i a_j\otimes v_i^j\in D_n^{\max}(V)$ for each $j$,
we have that $\sum_{j,i} a_j\otimes v_i^j\otimes w_i^j\in D_n$.
Thus $D_n^{\max}(V\otimes W)\subseteq D_n$.

For the reverse inclusion, the compatibility of the family
$\{D_n^{\max}(V\otimes W)\}_{n=1}^{\infty}$ implies that it suffices to show that if
$P\in D_k^{\max}(V)$ and $Q\in D_m^{\max}(W)$ then $P\otimes Q\in D_{km}^{\max}(V\otimes W)$.
However, such a $P$ (respectively, $Q$) has the form $P = \sum_{i=1}^l a_i\otimes v_i$
(respectively, $Q = \sum_{j=1}^r b_j\otimes w_j$), where $a_i\in M_k^+$ and $v_i\in V^+$
(respectively, $b_j\in M_m^+$ and $w_j\in W^+$), and hence
$$P\otimes Q = \sum_{i,j} (a_i\otimes b_j)\otimes (v_i\otimes w_j).$$
Clearly, $a_i\otimes b_j\in M_{km}^+$, and hence $D_n\subseteq D_n^{\max}(V\otimes W)$.
\endproof

\begin{remark}
{\rm
If $V$ and $W$ are AOU spaces, Effros defines in
\cite{Ef} their \lq\lq maximal tensor product'' $V \otimes_{MAX} W$ by using bilinear maps that are
  ``jointly positive''. (Effros actually uses lower case notation ``max'' for
  this tensor product, but we have adopted an upper case to avoid confusion.) Our jointly completely positive maps are the  ``complete'' analogue of these maps. In a recent preprint \cite{Han},
  Han also defines a maximal tensor product
  $V \otimes_{\pi} W$
  in the category of AOU spaces whose cone of positive elements coincides with our set
  $Q_{\max}.$ Combining \cite{Ef} with \cite{Han} (or just using
  \cite{Han}) one sees that these two definitions of the maximal
  tensor product in the category of AOU spaces coincide.  Thus Proposition \ref{p_omaxmax}
  shows that for any two AOU spaces $V$ and $W$ we have
$OMAX(V) \otimes_{\max} OMAX(W) = OMAX(V \otimes_{MAX} W).$ This
maximal tensor product of AOU spaces is also considered in Namioka
and Phelps \cite{np}.}
\end{remark}


\begin{remark}\label{th_ncsos}
{\rm Let A be a unital C*-algebra. Then $A$ is nuclear if and only if $A$ is $(\min, \max)$-nuclear;
that is, if and only if $A\otimes_{\min}\cl S = A\otimes_{\max}\cl S$ for every
operator system $\cl S$. Thus, the family of $(min,max)$-nuclear
operator systems contains the family of nuclear C*-algebras.}
\end{remark}
We give a proof here that relies on the Choi-Effors
characterization of nuclear C*-algebras, but in the next section we
will provide a proof that is independent of their result. The ``if'' part follows from Corollary \ref{c_minC} and
  Theorem \ref{th_maxcst}. To prove the converse implication we first
  show that $M_n\otimes_{\min} \cl S = M_n\otimes_{\max}\cl S$ for
  every operator system $\cl S$. In fact, we will show that these operator systems are
  both completely order isomorphic to $M_n(\cl S)$. If $\cl S$ is an operator
  subsystem of a C*-algebra $B$, then $M_n\otimes_{\min} \cl S$ is an operator
  subsystem of $M_n\otimes_{\min} B$ by injectivity. Note that
  $M_n\otimes_{\min} B = M_n(B)$, so $M_n\otimes_{\min} \cl S =
  M_n(\cl S)$. For the other equality note that if $u\in M_k(M_n(\cl
  S))^+$, then $u = \alpha (I_n\otimes u )\alpha^*$ where $\alpha =
  (E_{11}\,E_{21}\,\dots\,E_{n1})$ is in $M_k(M_n\otimes_{\max}\cl
  S)$. Since the cones of $M_n\otimes_{\max}\cl S$ are contained in those of $M_n\otimes_{\min}\cl S$, we obtain the
  desired equality.

  Now let $A$ be a nuclear C*-algebra. By
  \cite{CE3}, there exists a net of positive integers
  $\{n_{\lambda}\}$, unital completely positive maps $\phi_\lambda:A\rightarrow
  M_{n_\lambda}$, and unital completely positive maps $\psi_\lambda:M_{n_\lambda}
  \rightarrow A$ such that $\psi_\lambda \circ \phi_\lambda$
  converges to the identity on $A$ in the point-norm topology.

Consider the following maps:
$$
A\otimes_{\min} {\cl S} \xrightarrow{\phi_\lambda \otimes \id}  M_{n_\lambda} \otimes_{\min} {\cl S} \xrightarrow{\id_\lambda}  M_{n_\lambda} \otimes_{\max} {\cl S} \xrightarrow{\psi_\lambda \otimes \id} A\otimes_{\max} {\cl S},
$$
and let $\varphi_{\lambda}:A\otimes_{\min} {\cl S} \rightarrow
A\otimes_{\max} {\cl S}$ be their composition. More
precisely, $\nph_{\lambda}$ is given by $\varphi_{\lambda}(a\otimes s) = (\psi_{\lambda}\circ
\phi_{\lambda})(a)\otimes s$. Note that $\varphi_{\lambda}$ is
unital and completely
positive since the maps $\phi_\lambda \otimes \id$, $\id_\lambda$ and $\psi_\lambda \otimes \id$ are such.
We also observe that
$\varphi_{\lambda}$ approximates the identity in the sense that for every
$u\in A \otimes_{\max} \cl S$, we have
$\|\varphi_{\lambda}(u)-u\|\xrightarrow{\lambda} 0$. Indeed, if $u =
a\otimes s$ then  $\|\varphi_{\lambda}(a\otimes s) - a\otimes s\| =
\|(\psi_{\lambda}\circ \phi_{\lambda})(a)\otimes s - a\otimes s\| =
\|[(\psi_{\lambda}\circ \phi_{\lambda})(a) - a]\otimes s\| \leq
\|(\psi_{\lambda}\circ \phi_{\lambda})(a)-a\|\|
s\|\xrightarrow{\lambda} 0,$ where the inequality follows from the
fact that the operator space structure on $A \otimes_{\max} \cl S$ induces
an operator space cross-norm by Proposition~\ref{opsysopsp}. So the
result follows from the sublinearity of the norm.

Now let $U\in M_n( A\otimes_{\min} {\cl S})^+$.
Then $\varphi_{\lambda}^{(n)} (U) \in M_n( A\otimes_{\max} {\cl S})^+$ for every $\lambda$ and
$\varphi_{\lambda}^{(n)} (U)\rightarrow U$.
So we have that $M_n( A\otimes_{\min} {\cl S})^+ \subseteq M_n( A\otimes_{\max} {\cl S})^+$ since $M_n( A\otimes_{\max} {\cl S})^+$ is closed
by \cite[Theorem 2.30]{pt}.
The reverse inclusion is trivial.

We thus see that a C*-algebra is nuclear if and only if
$C_n^{\min}(A, \cl S) = C_n^{\max}(A, \cl S)$ for every $n \in \bb N$
and every operator system $\cl S.$

By Proposition \ref{th_ncsos}, every finite-dimensional C*-algebra
is ($\min$-$\max$)-nuclear. Unlike C*-algebras, finite-dimensional
operator systems do not have to be ($\min$-$\max$)-nuclear, as we
have observed in Remark~\ref{rem-fdopsys}. We now exhibit an
operator system that is \lq\lq nuclear'' when tensored with any
C*-algebra, but is not ($\min$, $\max$)-nuclear and is also not
(completely order isomorphic to) a C*-algebra. The operator system
defined in Theorem \ref{th_opse} will be fixed for the rest of this
section.

\begin{theorem}\label{th_opse}
Let $\cl S= span \{ E_{1,1}, E_{1,2}, E_{2,1}, E_{2,2},
  E_{2,3}, E_{3,2}, E_{3,3} \} \subseteq M_3.$  Then $\cl S
  \otimes_{\min} A = \cl S \otimes_{\max} A$ for every C*-algebra $A$, and $\cl S$ is not
  completely order isomorphic to a C*-algebra.
\end{theorem}
\begin{proof} By the injectivity of the minimal tensor product, we have that $\cl
  S \otimes_{\min} A \subseteq M_3 \otimes_{\min} A = M_3(A).$  Thus, to
  show that $C_n^{\max}(\cl S, A) = C_n^{\min}(\cl S, A),$ after
  identifying $M_n(\cl S \otimes A) = \cl S \otimes M_n(A),$ it will suffice to show that if
\[ P = \begin{pmatrix} P_{1,1} & P_{1,2} & 0\\ P_{2,1} & P_{2,2} &
  P_{2,3}\\ 0 & P_{3,2} & P_{3,3} \end{pmatrix} \in M_3(M_n(A))^+, \]
then $P \in C_n^{\max}.$

For every $r > 0$ we have that $r I_{n} +P_{i,i} > 0$ and that
\begin{multline*} rI_{3n} + P = \begin{pmatrix} rI_n+P_{1,1} & P_{1,2} & 0\\P_{2,1} &
  P_{2,1}(rI_n + P_{1,1})^{-1}P_{1,2} & 0\\0 & 0 & 0 \end{pmatrix}
+\\ \begin{pmatrix} 0 & 0 & 0\\0 & rI_n + P_{2,2} - P_{2,1}(rI_n +
  P_{1,1})^{-1}P_{1,2} & P_{2,3}\\0 & P_{3,2} & rI_n +
  P_{3,3} \end{pmatrix}. \end{multline*} Moreover, by the Cholesky algorithm both block matrices appearing in the sum
are positive.

By the nuclearity of $M_2$ and Theorem \ref{th_maxcst}, these matrices belong to
$C_n^{\max}(\cl S, A).$



To finish the proof we need to show that $\cl S$ is not completely
order isomorphic to a C*-algebra. Assume, by way of contradiction, that
$\cl S$ is completely order
isomorphic to a C*-algebra. Since $\dim(\cl S) =7,$ it must be completely order isomorphic to either $M_2 \oplus \bb C \oplus
\bb C \oplus \bb C$ or $\bb C \oplus \cdots \oplus \bb C$.
Since these C*-algebras are injective, $\cl S$ is injective.
This implies the existence of a completely
positive projection $\Psi$ from $M_3$ onto $\cl S.$
The map $\Psi$ fixes the algebra $\cl D_3$ of diagonal matrices and is hence a $\cl D_3$-bimodule map.
But such bimodule maps are given by Shur products with 3$\times$3 matrices.
It follows that $\Psi$ is given by Schur product against the
matrix
$R= \begin{pmatrix} 1 & 1 & 0\\1 & 1 & 1\\0 & 1 & 1 \end{pmatrix}.$
However, a Schur product map corresponding to a matrix $S$ is completely positive if and only if the
matrix $S$ is positive. Since $R$ is not a positive matrix, we obtain a contradiction which shows
that $\cl S$ can not be completely order isomorphic to a C*-algebra.
\end{proof}

We would like to point out that the fact that $\cl S$ is not
completely order isomorphic to a C*-algebra can also be deduced from
Theorem \ref{th_518} and Theorem \ref{th_maxcst}, but the above
argument avoids duality considerations.

We now wish to develop some further properties of the above operator system and of its dual. To this end, set \[G= \{(1,1), (1,2), (2,1), (2,2), (2,3), (3,2), (3,3) \}, \] so that $\cl S = \operatorname{span} \{ E_{i,j}: (i,j) \in G \}.$
Let $f_{i,j}: \cl S \to \bb C$, $i,j = 1,2,3$, be the dual functionals given by
$f_{i,j}(E_{k,l}) = \delta_{(i,j),(k,l)},$ where $\delta_{p,q}$ is the usual Kronecker delta function.
Then $\cl S^d = \operatorname{span} \{ f_{i,j}: (i,j) \in G \}.$

If $\cl T$ is an operator system and $f \in \cl T^d$ is a positive linear functional which is a matrix order unit for $\cl T^d$
it is easily seen that $f$ is Archimedean.
Thus, by \cite[Theorem 4.4]{CE2}, $(\cl T^d, \{ M_n(\cl T^d)^+ \}_{n=1}^\infty, f)$ is (completely order isomorphic to) an operator system.
It is shown in \cite[Corollary 4.5]{CE2} that whenever $\cl T$ is finite dimensional,
then such a functional $f$ exists and thus $\cl S^d$ is an operator system.
Below we give a concrete representation for $\cl S^d$.

\begin{proposition} Let $\cl S$ and $\cl S^d$ be as above, and let $A_{i,j} \in M_n$, $(i,j) \in G.$
Then $\sum_{(i,j) \in G} A_{i,j} \otimes f_{i,j} \in M_n(\cl S^d)^+$ if and only if $\begin{pmatrix} A_{1,1} & A_{1,2}\\A_{2,1} & A_{2,2} \end{pmatrix} \in M_2(M_n)^+$ and $\begin{pmatrix} A_{2,2} & A_{2,3}\\A_{3,2} & A_{3,3} \end{pmatrix} \in M_2(M_n)^+.$ Consequently, the linear map $\Gamma: \cl S^d \to M_2 \oplus M_2$ defined by
\[\Gamma\left(\sum_{(i,j) \in G} a_{i,j}f_{i,j}\right) =
\begin{pmatrix} a_{1,1}& a_{1,2}\\a_{2,1} & a_{2,2} \end{pmatrix} \oplus \begin{pmatrix} a_{2,2} & a_{2,3}\\a_{3,2} & a_{3,3} \end{pmatrix}, \]
is a complete order isomorphism onto its range.
\end{proposition}
\begin{proof} We have that
 $\sum_{(i,j) \in G} A_{i,j} \otimes f_{i,j}$ is in $M_n(\cl S^d)^+$ if and only if the map $\Phi: \cl S \to M_n$ defined by $\Phi(E_{i,j}) = A_{i,j}$ is completely positive.

If we assume that $\Phi$ is completely positive, then the restriction of $\Phi$ to $\operatorname{span} \{ E_{1,1}, E_{1,2}, E_{2,1}, E_{2,2} \} = M_2$ is completely positive. By a result of Choi, we have that $\begin{pmatrix} \Phi(E_{1,1}) & \Phi(E_{1,2}) \\ \Phi(E_{2,1}) & \Phi(E_{2,2}) \end{pmatrix} \in M_2(M_n)^+.$
In other words, $\begin{pmatrix} A_{1,1} & A_{1,2} \\ A_{2,1} & A_{2,2} \end{pmatrix} \in M_2(M_n)^+$.
Similarly, $\begin{pmatrix} A_{2,2} & A_{2,3} \\ A_{3,2} & A_{3,3} \end{pmatrix}$
can be seen to be positive by restricting to $\operatorname{span} \{ E_{2,2}, E_{2,3}, E_{3,2}, E_{3,3} \}.$

Conversely, if we assume that $\begin{pmatrix} A_{1,1} & A_{1,2} \\ A_{2,1} & A_{2,2} \end{pmatrix}$ and
$\begin{pmatrix} A_{2,2} & A_{2,3} \\ A_{3,2} & A_{3,3} \end{pmatrix}$ are positive,
then by the positive completion results of \cite{pr},
there exist $A_{1,3}, A_{3,1} \in M_n,$ such that $(A_{i,j})_{i,j=1}^3 \in M_3(M_n)^+.$ If we define $\Psi: M_3 \to M_n,$ via $\Psi(E_{i,j}) = A_{i,j},$ then we will have that $(\Psi(E_{i,j})) \in M_3(M_n)^+$ and so
again by Choi's result, $\Psi$ is completely positive.
Hence $\Phi$ is completely positive, since it is the restriction of $\Psi$ to an operator subsystem of $M_3$.
\end{proof}

\begin{theorem}\label{th_518}
The following hold for the operator system $\cl S$ and its dual $\cl S^d$:
\begin{enumerate}
\item If $A \subseteq B$ are unital C*-algebras and $\phi:A \to \cl S^d$ is completely positive, then $\phi$ possesses a completely positive extension $\psi:B \to \cl S^d.$
\item The identity map $\id: \Gamma(\cl S^d) \to \Gamma(\cl S^d)$ is a completely positive map that has no completely positive extension
to a map from $M_2 \oplus M_2$ to $\Gamma(\cl S^d).$
\item  $\id:\Gamma(\cl S^d) \otimes_{\min} \cl S \to \Gamma(\cl S^d) \otimes_{\max} \cl S$ is not completely positive.
\item $\cl S$ is not $(\min,\max)$-nuclear.
\end{enumerate}
\end{theorem}
\begin{proof}
By Theorem \ref{th_opse} and the fact that $\min$ and $\max$ are symmetric,
$A \otimes_{\max} \cl S = A \otimes_{\min}
\cl S \subseteq B \otimes_{\min} \cl S = B \otimes_{\max} \cl S,$
completely order isomorphically. Hence every jointly completely
positive map defined on $A \times \cl S$ can be extended to a jointly
completely positive map defined on $B \times \cl S.$ Part (1) now follows by
identifying $\phi:A \to \cl S^d$ with a jointly completely positive
map into $\bb{C}$, extending it to a jointly completely positive map from $B
\times \cl S$ into $\bb{C}$, and letting $\psi : B \to \cl S^d$ be the corresponding
linear map (see Lemma \ref{llance}).

To prove (2), suppose that the identity map on $\Gamma(\cl S^d)$ had a
completely positive extension $\Phi : M_2\oplus M_2\rightarrow \Gamma(\cl S^d)$.
Then $\Phi$ would be a completely positive projection onto $\Gamma(\cl S^d).$
We identify $M_2 \oplus M_2$ with the algebra of block diagonal matrices in $M_4$.
Under this identification, $\Gamma(f_{1,1}) = E_{1,1}$, $\Gamma(f_{1,2}) = E_{1,2}$,
$\Gamma(f_{2,1}) = E_{2,1}$, $\Gamma(f_{2,2}) = E_{2,2} + E_{3,3}$,
$\Gamma(f_{2,3}) = E_{3,4}$, $\Gamma(f_{3,2}) = E_{4,3}$, and
$\Gamma(f_{3,3}) = E_{4,4}$. Thus, $\cl D =
\operatorname{span} \{E_{1,1}, E_{2,2}+E_{3,3}, E_{4,4} \}$ would be a C*-algebra fixed by
$\Phi$, and hence $\Phi$ would be a $\cl D$-bimodule map (see \cite[Corollary 3.19]{Pa}).
Since $\Phi(E_{2,2}) \in \Gamma(\cl S^d)$
and $(E_{2,2}+E_{3,3})\Phi(E_{2,2}) = \Phi(E_{2,2}) = \Phi(E_{2,2})(E_{2,2}+E_{3,3}),$ we would have that
$\Phi(E_{2,2}) = t(E_{2,2} +E_{3,3})$ for some $t \ge 0.$  Similarly, $\Phi(E_{3,3}) = r(E_{2,2} +E_{3,3})$
for some $r \ge 0,$ and it would follow that $t+r=1.$  But since
$0 \le J_1 = E_{1,1} +E_{1,2} + E_{2,1} + E_{2,2}$, we have that
$0 \le \Phi(J_1) = E_{1,1} + E_{1,2} + E_{2,1} + tE_{2,2},$ and hence $t=1.$
Similarly, considering $J_2= E_{3,3} + E_{3,4} +E_{4,3} +E_{4,4}$ yields
that $r=1,$ contradicting the fact that $r+t=1.$

To see (3), suppose that the identity map is completely positive.  Then we have that $\Gamma(\cl S^d) \otimes_{\max} \cl S = \Gamma(\cl S^d) \otimes_{\min}  \cl S \subseteq (M_2 \oplus M_2) \otimes_{\min} \cl S = (M_2 \oplus M_2) \otimes_{\max} \cl S,$ where the identifications and inclusions are in the complete order sense. These inclusions imply that every jointly
completely positive map on $\Gamma(\cl S^d) \times \cl S$ extends to a jointly
completely positive map on $(M_2 \oplus M_2) \times \cl S.$  Thus every completely positive map from $\Gamma(\cl S^d)$ into $\cl S^d = \Gamma(\cl S^d)$ extends to a completely positive map from $M_2 \oplus M_2$ to $\Gamma(\cl S^d)$, which contradicts (3).

(4) is a direct consequence of (3).
\end{proof}

The above results show that even though $A \otimes_{\min} \cl S = A \otimes_{\max} \cl S$ for every C*-algebra, neither $\cl S$ nor $\cl S^d$ is injective.

\begin{remark}\label{graph_rem}
{\rm A graph $G$ on $n$ vertices can be identified with a subset
$G \subseteq \{1,\dots,n \} \times \{ 1,\dots, n \}$ satisfying the properties that $(i,j) \in
G$ whenever $(j,i) \in G$ and that $(i,i) \in G$ for $i=1,\dots,n.$
To such a graph one can associate an operator system $\cl S(G) = \operatorname{span}
\{ E_{i,j} : (i,j) \in G \} \subseteq M_n.$  One can show
that if the graph G is chordal, then $\cl S(G) \otimes_{\min} A = \cl
S(G) \otimes_{\max} A$ for every C*-algebra $A.$ The proof is similar
to that of Theorem \ref{th_opse} and uses the fact that chordal graphs have a
``perfect vertex elimination scheme'' and the techniques of
\cite{pps} and \cite{pr},
where it is shown that
whenever one has a perfect vertex elimination scheme, then one can
carry out a Cholesky-type algorithm as above to decompose strictly
positive matrices in $\cl S(G) \otimes_{\min} A$ as encountered in the
proof of Theorem \ref{th_opse}. We do not present this argument here
though, since this result also follows more readily from results in
the next section.

We note that the operator system $\cl S$ of Theorem~\ref{th_opse} is the operator system associated to the following chordal graph:
$$\xymatrix{ \bullet \ar@{-}[r] & \bullet \ar@{-}[r] & \bullet
}$$
}
\end{remark}


\section{The commuting tensor product}

In this section we introduce another operator system tensor product which
agrees with the $\max$ tensor product for all pairs of C*-algebras,
but does not agree with the $\max$ tensor product on all pairs of operator systems.
Thus, this new operator system tensor product gives a different extension of the
maximal C*-algebraic tensor product from the category of C*-algebras to the category of
operator systems.  In contrast with the maximal operator system tensor product,
but in analogy with the minimal one, this tensor product is defined by specifying
a collection of completely positive maps rather than specifying the matrix ordering.

Let $\cl S$ and $\cl T$ be operator systems. Set
\begin{align*}
\cp(\cl S,\cl T) = \{(\phi,\psi) : &\text{ $H$ is a Hilbert space,  $\phi\in \operatorname{CP}(\cl S, \cl B(H))$,} \\
& \text{$\psi \in \operatorname{CP}(\cl T, \cl B(H))$, and $\phi(\cl S)$ commutes with $\psi(\cl T)$.} \}
\end{align*}
Given $(\phi,\psi)\in \cp(\cl S,\cl T)$, let $\phi\cdot\psi : \cl S\otimes\cl T\rightarrow \cl B(H)$
be the map given on elementary tensors by $(\phi\cdot\psi) (x\otimes y) = \phi(x)\psi(y)$.

For each $n\in \bb{N}$, define a cone $P_n \subseteq M_n(\cl S\otimes\cl T)$ by letting
$$P_n = \{u\in M_n(\cl S\otimes \cl T) : (\phi\cdot\psi)^{(n)}(u) \geq 0, \mbox{ for all }
(\phi,\psi)\in \cp(\cl S,\cl T)\}.$$

\begin{proposition}\label{p_c}
The collection $\{P_n \}_{n=1}^{\infty}$ is a matrix ordering on $\cl S\otimes\cl T$
with Archimedean matrix unit $1\otimes 1$.
\end{proposition}

\proof It is clear that $P_n$ is a cone. If $\alpha\in M_{n,m}$ and $u\in P_m$ then
$$(\phi\cdot\psi)^{(n)}(\alpha u \alpha^*) = \alpha (\phi\cdot\psi)^{(m)}(u)\alpha^* \geq 0,$$
and hence the family $\{P_n \}_{n=1}^{\infty}$ is compatible.
Let $\phi\in S_k(\cl S)$ and $\psi\in S_m(\cl T)$, and define
$\tilde{\phi} : \cl S\rightarrow M_k\otimes 1_m$
(respectively, $\tilde{\psi} : \cl T\rightarrow 1_k\otimes M_m$) by
$\tilde{\phi}(x) = \phi(x)\otimes 1_m$ (respectively, $\tilde{\psi}(y) = 1_k \otimes \psi(y)$).
Then $(\tilde{\phi},\tilde{\psi})\in \cp(\cl S,\cl T)$ and hence
$$(\phi\otimes\psi)^{(n)}(u) = (\tilde{\phi}\cdot\tilde{\psi})^{(n)}(u)\geq 0 \quad \text{ for each }  u\in P_n.$$
Thus $P_n\subseteq C_n^{\min}$ for each $n\in \bb{N}$.
It now follows that $P_n \cap (-P_n) = \{0\}$ and that $1\otimes 1$ is an
matrix order unit for $\{P_n \}_{n=1}^{\infty}$.

Suppose that $r(1\otimes 1)_n + u\in P_n$ for each $r > 0$.
Then $(\phi\cdot\psi)^{(n)}(r(1\otimes 1)_n + u) \geq 0$ for all $(\phi,\psi)\in \cp(\cl S,\cl T)$ and all $r > 0$.
Thus $r I_{H} + (\phi\cdot\psi)^{(n)}(u) \geq 0$ for all $(\phi,\psi)\in \cp(\cl S,\cl T)$ and all $r > 0$,
which implies that $u\in P_n$. Hence, $1\otimes 1$ is an Archimediean matrix order unit, and the proof is complete.
\endproof

\begin{definition}\label{d_c}
We let $\cl S\otimes_{\comm}\cl T$ denote the operator system $(\cl S\otimes\cl T, \{ P_n\}_{n=1}^{\infty}, 1 \otimes 1)$.
\end{definition}

\begin{theorem}\label{th_c}
The mapping $\comm : \cl O\times\cl O\rightarrow\cl O$ sending the
pair $(\cl S,\cl T)$ to the operator system $\cl S\otimes_{\comm}\cl
T$ is a symmetric, functorial operator system tensor product.
\end{theorem}
\proof
Properties~(T1) and (T3) were checked in the proof of Proposition \ref{p_c}.
Suppose $P = (p_{i,j})\in M_n(\cl S)^+$ and $Q \in M_m(\cl T)^+$, and let
$(\phi,\psi)\in \cp(\cl S,\cl T)$. Then
\begin{eqnarray*}
(\phi\cdot\psi)^{(nm)}(P\otimes Q) & = & ((\phi\cdot\psi)^{(m)}(p_{i,j}\otimes Q))_{i,j}\\
& = & ((\phi(p_{i,j})\otimes 1_m)\psi^{(m)}(Q))_{i,j}\\ & = &
(\phi(p_{i,j})\otimes 1_m)_{i,j}(\psi^{(m)}(Q)\otimes 1_n) \geq 0.
\end{eqnarray*}
It follows that Property~(T2) is satisfied, and hence $\comm$ is an operator system tensor product.

We next check functoriality. Suppose that $\rho : \cl S_1\rightarrow \cl S_2$ and
$\eta : \cl T_1\rightarrow\cl T_2$ are unital completely positive maps,
and let $u\in M_n(\cl S_1\otimes_{\comm}\cl T_1)^+$. If $(\phi',\psi')\in \cp(\cl S_2,\cl T_2)$, then $(\phi'\circ \rho,\psi'\circ\eta)\in \cp(\cl S_1,\cl T_1)$, and
$$(\phi'\cdot\psi')^{(n)}((\rho\otimes\eta)^{(n)}(u)) = ((\phi'\circ \rho)\cdot(\psi'\circ\eta))^{(n)}(u) \geq 0,$$
and hence $(\rho\cdot\eta)^{(n)}(u)$ is in the positive cone of $M_n(\cl S_2\otimes_{\comm}\cl T_2)$.
This establishes Property~(T4).

Recall that $\theta : \cl S\otimes\cl T\rightarrow\cl T\otimes\cl S$ denotes the
map given by $\theta(x\otimes y) = y\otimes x$.
Note that $(\phi,\psi)\in \cp(\cl S,\cl T)$ if and only if $(\psi,\phi)\in \cp(\cl T,\cl S)$.
Moreover, if $u\in \cl S\otimes\cl T$
and $(\phi,\psi)\in \cp(\cl S,\cl T)$, then
$(\phi\cdot\psi)(u) = (\psi\cdot\phi)(\theta(u))$. It follows
that the tensor product $\comm$ is symmetric.
\endproof

We recall that for an operator system $\cl S$, there exists a
unital C*-algebra $C^*_u(\cl S)$ introduced in \cite{KW} (called either the
\textbf{universal C*-algebra of $\cl S$} or the \textbf{maximal C*-algebra of $\cl S$})
and a unital completely positive map $\iota : \cl
S\rightarrow C_u^*(\cl S)$ with the properties that $\iota(\cl S)$
generates $C^*_u(\cl S)$ as a C*-algebra, and that for every unital
completely positive map $\phi : \cl S\rightarrow \cl B(H)$ there
exists a unique $*$-homomorphism $\pi : C^*_u(\cl S) \to \cl B(H)$ such that
$\pi \circ \iota = \phi$. To construct this C*-algebra, one starts with the free $*$-algebra
$$\cl F(\cl S) = \cl S \oplus (\cl S \otimes \cl S) \oplus
(\cl S \otimes \cl S \otimes \cl S) \oplus \cdots.$$ Each unital
completely positive map $\phi: \cl S \to \cl B(H)$ gives rise to a
$*$-homomorphism $\pi_{\phi}: \cl F(\cl S) \to \cl B(H)$ by setting
\[\pi_{\phi}(s_1 \otimes \cdots \otimes s_n) = \phi(s_1) \cdots \phi(s_n)\]
and extending linearly to the tensor product and then to the direct sum.
For $u\in \cl F(\cl S)$, one sets $\|u\|_{\cl F(\cl S)} = \sup\|\pi_{\phi}(u)\|$,
where the supremum is taken over all unital completely positive maps $\phi$ as above,
and defines $C^*_u(\cl S)$ to be the completion of $\cl F(\cl S)$ with respect to $\|\cdot\|_{\cl F(\cl S)}$.
We will identify $\cl S$ with its image $\iota(\cl S)$, and thus consider
$\cl S$ as an operator subsystem of $C^*_u(\cl S)$.


\begin{theorem}\label{th_cinu}
Let $\cl S$ and $\cl T$ be operator systems. The operator system
arising from the inclusion of $\cl S\otimes\cl T$ into $C^*_u(\cl
S)\otimes_{\max} C^*_u(\cl T)$ coincides with $\cl
S\otimes_{\comm}\cl T$.
\end{theorem}
\proof
Let $\tau$ be the operator system structure on $\cl S\otimes\cl T$
arising from the inclusion $\cl S \otimes\cl T\subseteq C^*_u(\cl S)\otimes_{\max} C^*_u(\cl T)$.
Suppose that $u\in M_n(\cl S\otimes_{\tau}\cl T)^+$ and let $(\phi,\psi)\in \cp(\cl S,\cl T)$.
By the universal properties of $C^*_u(\cl S)$ and $C^*_u(\cl T)$, there exist (unique)
$*$-homomorphisms $\pi : C^*_u(\cl S)\rightarrow \cl B(H)$ and $\rho : C^*_u(\cl T)\rightarrow \cl B(H)$
extending $\phi$ and $\psi$, respectively. Since $\cl S$ (respectively, $\cl T$) generates
$C^*_u(\cl S)$ (respectively, $C^*_u(\cl T)$) as a C*-algebra, we have that the ranges of
$\pi$ and $\rho$ commute. It follows that
$$(\phi\cdot\psi)^{(n)}(u) = (\pi\cdot\rho)^{(n)}(u)\geq 0,$$
and hence $u\in M_n(\cl S\otimes_{\comm}\cl T)$.

Conversely, suppose that $u\in M_n(\cl S\otimes_{\comm}\cl T)^+$. To
show that $u$ is in the positive cone of $M_n(C^*_u(\cl S)\otimes_{\max}C^*_u(\cl T))$, it
suffices by Lemma \ref{l_gotoa} to prove that $\eta^{(n)}(u)\geq 0$ for each
completely positive map $\eta : C^*_u(\cl S)\otimes_{\max} C^*_u(\cl T)\rightarrow \cl
B(H)$. By Stinespring's Theorem, we may moreover assume that $\eta$ is a
$*$-homomorphism.
By Theorem~\ref{th_maxcst} and the universal
property of the maximal tensor product of C*-algebras, each such $\eta$
is equal to $\pi\cdot\rho$, where $\pi : C^*_u(\cl
S)\rightarrow \cl B(H)$ and $\rho : C^*_u(\cl T)\rightarrow\cl
B(H)$ are $*$-homomorphisms with commuting ranges. Since the
restrictions of $\pi$ to $\cl S$ and the restriction of $\rho$ to $\cl T$ are each completely positive, we have that $\eta(u)\geq 0$.
\endproof

We obtain the following consequence of Theorem \ref{th_cinu}.


\begin{corollary}\label{c_chcccom}
Let $\cl S$ and $\cl T$ be operator systems. A linear map $f : \cl
S\otimes_{\comm}\cl T \rightarrow \cl B(H)$ is a unital completely
positive map if and only if there exist a Hilbert space $K$,
$*$-homomorphisms $\pi : C^*_u(\cl S)\rightarrow \cl B(K)$ and
$\rho : C^*_u(\cl T)\rightarrow \cl B(K)$ with commuting ranges, and
an isometry $V : H\rightarrow K$ such that $f(x \otimes y) =
V^*\pi(x)\rho(y)V$ for all $x\in \cl S$ and all $y\in \cl T$.
\end{corollary}
\proof Suppose that $K$, $V$, $\pi$, and $\rho$ are as in the
statement. Since $\pi\cdot\rho$ is completely positive on $C^*_u(\cl
S)\otimes_{\max} C^*_u(\cl T)$, Theorem~\ref{th_cinu} implies that
the restriction of $\pi\cdot\rho$ to $\cl S\otimes_{\comm}\cl T$ is completely
positive. Hence the map $u\rightarrow V^*\pi\cdot\rho(u)V$ on $\cl
S\otimes_{\comm}\cl T$ is completely positive.

Conversely, suppose that $f : \cl S\otimes_{\comm}\cl T\rightarrow
\cl B(H)$ is completely positive. By Theorem~\ref{th_cinu}, $f$
has a completely positive extension $\tilde{f} : C^*_u(\cl
S)\otimes_{\max} C^*_u(\cl T)\rightarrow \cl B(H)$. Stinespring's
Theorem implies the existence of a Hilbert space $K$, an isometry $V
: H\rightarrow K$, and a *-homomorphism $\eta : C^*_u(\cl
S)\otimes_{\max} C^*_u(\cl T)\rightarrow \cl B(K)$ such that
$\tilde{f}(u) = V^*\eta(u)V$ for all $u\in C^*_u(\cl S)\otimes_{\max}
C^*_u(\cl T)$. By the universal property of the maximal C*-algebraic
tensor product, $\eta = \pi\cdot\rho$ for some $*$-homomorphisms $\pi
: C^*_u(\cl S)\rightarrow \cl B(K)$ and $\rho : C^*_u(\cl
T)\rightarrow \cl B(K)$.
\endproof

The next result, Theorem~\ref{th_coincm}, can be deduced as a
corollary of the following Theorem~\ref{stronger-thm}, but we
present a separate proof because it is a considerably more elementary result.

\begin{theorem}\label{th_coincm}
If $A$ and $B$ are unital C*-algebras, then $A\otimes_{\comm} B = A\otimes_{\max} B$.
\end{theorem}
\begin{proof} By Theorem~\ref{th_max}, $M_n(A\otimes_{\max} B)^+\subseteq
M_n(A\otimes_{\comm} B)^+$. Conversely, suppose that $u\in
M_n(A\otimes_{\comm} B)^+$. By Theorem \ref{th_maxcst},
$A\otimes_{\max} B$ is completely order isomorphic to the image of
$A\otimes B$ inside $A\otimes_\text{C*max}B$, the maximal
C*-algebraic tensor product of $A$ and $B$. Now let
$i_A:A\rightarrow A\otimes_\text{C*max}B$ be given by $i_A (a) =
a\otimes {1_B}$ and let $i_B:B\rightarrow A\otimes_\text{C*max}B$ be
given by $i_B(b) = 1_A\otimes b$. Clearly, $i_A$ and $i_B$ are
completely positive and have commuting ranges.
Theorem~\ref{th_maxcst} implies that $u\in M_n(A\otimes_{\max} B)^+$
if and only if $(i_A\cdot i_B)^{n}(u)$ is positive. But the latter
is true by the definition of the commuting tensor product. Thus the
result follows.
\end{proof}

The following result gives another characterization of the $\comm$
tensor product. We prove that, in a certain precise sense, $\comm$ is
the minimal extension of $C^*max$ from the category of C*-algebras to
the category of operator systems.

\begin{theorem} \label{stronger-thm}
If $A$ is a unital $C^*$-algebra and $\cl S$ is an operator system,
then $A \otimes_{\comm} \cl S=A \otimes_{\max} \cl S$. Moreover, if
$\alpha: \cl O \times \cl O \to \cl O$ is any symmetric, functorial
operator system tensor product such that $A \otimes_{\alpha} B = A
\otimes_{\max} B$ for every pair of unital $C^*$-algebras $A,$ then
$\comm \le \alpha,$ i.e., for every pair of operator systems $\cl S$
and $\cl T,$ the
identity map $id_{\cl S} \otimes id_{\cl T}: \cl S \otimes_{\alpha} \cl T
\to \cl S \otimes_{\comm} \cl T$ is completely positive.
\end{theorem}
\begin{proof}
By defining $a_1 \cdot (a \otimes s) \cdot a_2 = (a_1aa_2) \otimes s$,
the algebraic tensor product $ A \otimes \cl S$ becomes an $
A$-bimodule. We claim that $ A \otimes_{\max} \cl S$ is an {\it
  operator $ A$-system} in the sense of \cite[Chapter~15]{Pa}; that
is, if $U\in M_n( A \otimes_{\max} \cl S)^+$ and $B\in M_{n,k}( A)$,
then $B^*\cdot U \cdot B$ is in $M_k( A \otimes_{\max} \cl S)^+$. To
show this, we may assume that $U$ is in $D_n^{\max}$.  Indeed, suppose that the
assertion is true in this case. Given $V \in C_n^{\max}$, we know
that $V+\epsilon I_n \in D_n^{\max}$ for every $\epsilon >0$. We have that
$B^*\cdot (V+\epsilon I_n) \cdot B = B^*\cdot V \cdot B + \epsilon
B^*\cdot I_n \cdot B = B^*\cdot V \cdot B + \epsilon B^*B\otimes
(1_{\cl S})$ is in $C_n^{\max}$ for every $\epsilon >0$. So the result
follows from the fact that $C_n^{\max}$ is closed.

Let $U\in D_n^{\max}$ have the form
$U = \alpha (P\otimes Q) \alpha^*$,
where $P \in M_p( A)^+$, $Q = (s_{ij})\in M_q(\cl S)^+$ and $\alpha \in M_{n,pq}$.
Note that $B^*\cdot \alpha (P\otimes Q) \alpha^*\cdot B = (\alpha^* B)^*\cdot (P\otimes Q)\cdot (\alpha^* B)$.
Thus we may assume that $U = P\otimes Q$, where $P \in M_p( A)^+$ and $Q = (s_{ij})\in M_q(\cl S)^+$ with $pq=n$.
Let $B = (B_1\,B_2 \ldots B_q)^{\ttt}$, where each $B_i$ is a $p\times k$ matrix. Then
\begin{align*}
B^*\cdot(P\otimes Q)\cdot B &= (B_1^*\,B_2^*\,...\,B_q^*)\cdot
\left(
\begin{array}{ccc}
P\otimes s_{11} & \cdots & P\otimes s_{1q}\\
\vdots          & \ddots & \vdots   \\
P\otimes s_{q1} & \cdots & P\otimes s_{qq}
\end{array}
\right)\cdot
\left(
\begin{array}{c}
B_1    \\
\vdots \\
B_q
\end{array}
\right) \\
&= \sum_{i,j=1}^q (B_i^* P B_j) \otimes s_{ij}.
\end{align*}

Let $C = (B_i^* P B_j)_{i,j=1}^q,$ so that $C\in M_{kq}( A)^+$.  Let $X = (e_1\otimes I_k\,...\,e_q\otimes I_k)^{\ttt}$, where $e_i\otimes I_k = (0\,\dots\,I_k\,\dots\,0)^{\ttt}$ is a $qk\times k$ scalar matrix.  Then
$$
B^*\cdot(P\otimes Q)\cdot B = X^* (C\otimes Q)X \in C_{n}^{\max}.
$$

Thus we have shown that $A \otimes_{\max} \cl S$ is an operator $
A$-system. By \cite[Theorem~15.12]{Pa}, the map
$\pi: A \rightarrow \cl I( A \otimes_{\max} \cl S)$
given by $\pi(a) = a\otimes 1_{\cl S}$ is a unital
$*$-homomorphism. In this case, $\pi$ is also injective and hence an isometry.

Let $i:\cl S \rightarrow \cl I( A \otimes_{\max} \cl S)$ be given by
$i(s) = 1_{A}\otimes s$.  Then $i$ is a complete order isomorphism
onto its range. Note that $\pi( A)$ commutes with
$i(\cl S)$ since $ (a \otimes 1_{\cl S})(1_{ A}\otimes s) = a\cdot
(1_{A} \otimes 1_{\cl S})(1_{ A}\otimes s) =a\cdot (1_{ A}\otimes s) =
a\otimes s =  (1_{ A}\otimes s)\cdot a = (1_{ A}\otimes s)(1_{ A}
\otimes 1_{\cl S})\cdot a = (1_{ A}\otimes s)(a \otimes 1_{\cl S})
$. Thus $\pi: A\rightarrow \cl I( A \otimes_{\max} \cl S)$ and $i:\cl S
\rightarrow \cl I( A \otimes_{\max} \cl S)$ are completely
positive and have
commuting ranges. This means that $\pi\cdot i:  A \otimes_{\comm} \cl
S \rightarrow \cl I( A \otimes_{\max} \cl S)$ is completely
positive with range $
A \otimes_{\max} \cl S$. Note that $\pi\cdot i (a\otimes s) = a\otimes
s$, which implies that the identity map from $ A \otimes_{\comm} \cl S$
to $ A \otimes_{\max} \cl S$ is completely
positive. Thus $ A \otimes_{\comm} \cl S =
A \otimes_{\max} \cl S$ by the maximality of $\max$.

Finally, to see the last claim assume that $\alpha$ is as above and
let $\cl S$ and $\cl T$ be operator systems. By the functoriality of
$\alpha$ we have that the inclusion maps $\cl S \to C^*_u(\cl S)$ and
$\cl T \to C^*_u(\cl T)$ induce a completely positive map
$\cl S \otimes_{\alpha} \cl T \to C^*_u(\cl S) \otimes_{\alpha} C^*_u(\cl
T) =  C^*_u(\cl S) \otimes_{max} C^*_u(\cl T).$ But we also have that
these inclusion maps induce a complete order isomorphism of $\cl S
\otimes_{\comm} \cl T$ into $C^*_u(\cl S) \otimes_{\max} C^*_u(\cl T)$
and the result follows.
\end{proof}

We can now give the promised proof of Remark~5.15.

\begin{corollary} Let A be a unital $C^*$-algebra. Then A is a nuclear $C^*$-algebra if and
  only if $A \otimes_{\min} \cl S = A \otimes_{\max} \cl S$ for every
  operator system $\cl S.$
\end{corollary}
\begin{proof} We have that $A \otimes_{\max} \cl S = A \otimes_{\comm}
  \cl S \subseteq A \otimes_{\max} C^*_u(\cl S) = A \otimes_{\min}
  C^*_u(\cl S)$ and $A \otimes_{\min} \cl S \subseteq A \otimes_{\min}
  C^*_u(\cl S),$ where both containments are complete order
  isomorphisms. Thus, the result follows.
\end{proof}

We now define a tensor product for operator spaces that is related to
the $\mu$ tensor product of Oikhberg and Pisier\cite{op}. Let $X$ and $Y$
be operator spaces. For $u\in X\otimes Y$, let
\begin{align*}
\|u\|_{\mu^*} = \sup\{\|(f&\cdot g)(u)\| \ : \ \text{$f : X\rightarrow \cl
B(H)$ and $g : Y\rightarrow \cl B(H)$ are} \\
&\text{completely contractive maps with the property that $f(x)$} \\
& \quad  \text{commutes with $\{g(y),g(y)^*\}$ for all $x\in X$ and $y\in Y$} \}.
\end{align*}
We define norms on $M_n(X\otimes Y)$ in a similar fashion.
It is easily checked that this gives an operator space structure to $X \otimes Y$, and we denote the resulting
operator space $X\otimes_{\mu^*}Y$.
If the mappings $f$ and $g$ satisfy the properties in the definition of $\|\cdot\|_{\mu^*}$, we say that
their ranges are $*$-commuting.

\begin{proposition}\label{p_resmst}
Let $X$ and $Y$ be operator spaces. Then the identity map is a completely isometric
isomorphism betwen $X\otimes^{\comm} Y$ and $X\otimes_{\mu^*} Y$.
\end{proposition}
\begin{proof}
Given unital completely positive maps $\Phi: \cl S_X \to \cl B(K)$
and $\Psi: \cl S_Y \to \cl B(K)$ with commuting ranges, define $f:X
\to \cl B(K)$ and $g:Y \to \cl B(K)$ via $f(x) = \Phi(
\begin{pmatrix} 0 &
x\\0 & 0 \end{pmatrix} )$ and $g(y) = \Psi( \begin{pmatrix} 0 & y\\
0 & 0
\end{pmatrix} ).$  Then $f$ and $g$ are completely contractive maps
whose ranges are $*$-commuting.  This shows that the norm on $X
\otimes_{\mu^*} Y$ is greater than the norm on $X \otimes^{\comm}
Y.$

Conversely, given completely contractive commuting maps $f:X \to \cl
B(H)$ and $g:Y \to \cl B(H)$ as in the above definition define
completely positive maps $\Phi: \cl S_X \to \cl B(H \oplus H \oplus
H \oplus H)$ and $\Psi: \cl S_Y \to \cl B(H \oplus H \oplus H \oplus
H)$ by
\[ \Phi( \begin{pmatrix} \lambda & x_1\\ x_2^* & \mu \end{pmatrix} ) =
\begin{pmatrix} \lambda I_H & f(x_1) & 0 & 0\\ f(x_2)^* & \mu I_H & 0 & 0 \\ 0 & 0 & \lambda I_H & f(x_1)\\0 & 0 & f(x_2)^* & \mu I_H \end{pmatrix} \]
and
\[ \Psi( \begin{pmatrix} \alpha  & y_1\\y_2^* & \beta \end{pmatrix}) =
\begin{pmatrix} \alpha I_H & 0 & g(y_1) & 0\\ 0 & \alpha I_H & 0 & g(y_1)\\ g(y_2)^* & 0 & \beta I_H & 0\\ 0 & g(y_2)^* & 0 & \beta I_H \end{pmatrix}.
\]
The maps $\Phi$ and $\Psi$ are readily seen to be unital completely positive and
to have commuting ranges. This shows that the norm on $X \otimes_{\mu^*} Y$ does not exceed
the norm on $X \otimes^{\comm} Y$, and hence the two norms are equal.
\end{proof}

\begin{corollary}\label{c_diffe}
The operator system tensor products $\max$ and $\comm$ are distinct.
\end{corollary}
\begin{proof}
It will be enough to show that the induced operator space tensor
products $\otimes^{\max}$ and $\otimes^{\comm}$ are different.  In
\cite{op} Oikhberg and Pisier  introduce a tensor norm
$\otimes_{\mu}$ on operator spaces by considering the supremum over
all pairs of commuting (but not necessarily $*$-commuting)
completely contractive maps, and prove that this tensor norm is
strictly smaller than the projective operator space tensor norm.
Clearly, our $\|\cdot\|_{\mu^*} = \|\cdot\|_{\comm}$ is dominated by
$\|\cdot\|_{\mu}$ and since, by Theorem \ref{th_resmaxi},
$\|\cdot\|_{\max}$ coincides with the operator projective tensor
norm, the result follows.
\end{proof}

For the next result we need to recall the operator systems associated
with graphs that were introduced in Remark~\ref{graph_rem}

\begin{proposition}\label{p_digraph}
Let $G \subseteq \{ 1,\dots,k \} \times \{ 1,\dots,k \}$ be a graph on $k$
vertices and let $\cl S(G) \subseteq M_k$ be the operator system of
the graph. If $G$ is a chordal graph, then $\cl S(G) \otimes_{\comm}\cl T = \cl S(G) \otimes_{\min}\cl T$ for every operator
system $\cl T,$ and so $\cl S(G)$ is $(\min,\comm)$-nuclear.
\end{proposition}
\proof
Let $\{E_{i,j}\}$
be the canonical matrix units in $M_k$. Suppose that $\phi :
\cl S(G) \rightarrow \cl B(H)$ and $\psi : \cl T\rightarrow \cl B(H)$
are completely positive maps with commuting ranges. Let $T_{i,j} =
\phi(E_{i,j})$, $(i,j) \in G$. For every complete subgraph $G_0\subseteq
G$ (that is, a subset $G_0$ of $G$ of the form
$G_0 = J\times J$ for some $J\subseteq \{1,\dots,k\}$), we have
that $\phi |_{\cl S(G_0)} : \cl S(G_0)\rightarrow \cl
B(H)$ is completely positive. It follows by Choi's characterization
\cite{Ch} that the matrix $(T_{i,j})_{(i,j)\in G_0}$ is positive.

Thus, the partially defined matrix $(T_{i,j})_{(i,j) \in G}$ is
{\it partially positive} in the sense of \cite{pr}. It follows from
\cite{pr} that this operator matrix has a positive
completion in the von Neumann algebra $\phi(\cl S(G))^{\prime
  \prime}$; that is, there exist $T_{i,j}\in
\phi(\cl S(G))''$ for $(i,j)\not\in G$, such that the (fully defined)
operator matrix
$(T_{i,j})_{i,j=1}^k$ is positive. Another application of Choi's Theorem implies that the mapping
$\tilde{\phi} : M_k\rightarrow \cl B(H)$ sending a matrix $(\lambda_{i,j})$
to the operator $\sum_{i,j=1}^k \lambda_{i,j}T_{i,j}$ is completely positive.
Thus, $\tilde{\phi}$ is a completely positive extension of $\phi.$
Clearly the ranges of $\tilde{\phi}$ and $\psi$ commute.

It follows from the previous paragraph that $\cl
S(G)\otimes_{\comm}\cl T\subseteq M_k\otimes_{\comm}\cl T$ as
operator systems. However, $M_k$ is a nuclear C*-algebra, and hence
Theorem~\ref{th_ncsos} implies that $M_k\otimes_{\comm}\cl T =
M_k\otimes_{\min}\cl T$. On the other hand, $\cl
S(G)\otimes_{\min} \cl T\subseteq M_k\otimes_{\min}\cl T$ by
the injectivity of the minimal operator system tensor product. It
follows that $\cl S(G)\otimes_{\comm}\cl T = \cl
S(G)\otimes_{\min}\cl T$.
\endproof

Combining this proposition with Theorem~\ref{stronger-thm}, we have
that when $G$ is a chordal graph and $A$ is a C*-algebra, then
\[ \cl S(G) \otimes_{\min} A = \cl S(G) \otimes_{\comm} A = \cl S(G)
\otimes_{\max} A, \]
which is the result claimed in Remark~\ref{graph_rem}.

It follows from Proposition \ref{p_digraph} that the 7 dimensional
operator system of Theorem \ref{th_518} is $(\min,\comm)$-nuclear but
not $(\min,\max)$-nuclear.


\section{The lattice of  tensor products}\label{s_incl}

In this section we examine the collection of all operator system
tensor products, show that it is a lattice, and introduce some tensor
products that can also be characterized via this lattice. These tensor
products appear to have important categorical roles and are natural
analogues of some of the tensor products that appear in Grothendieck's
programme. We then relate preservation of these tensor products to
certain important properties of C*-algebras. First we
will need a preliminary result.

\begin{proposition}\label{p_latt}
The collection of all operator system tensor products is a complete
lattice with respect to the order introduced in Section~\ref{s-gs}. The collection
of all functorial operator system tensor products is a complete
sublattice of this lattice.
\end{proposition}
\proof Let $\{\tau_j\}_{j\in J}$ be a collection of operator system
tensor products, where $J$ is a non-empty set. It suffices to show
that $\{\tau_j\}_{j\in J}$ possesses a greatest lower bound. Fix
operator systems $\cl S$ and $\cl T$. For each $n\in \bb{N}$, let
$P_n = \bigcap_{j\in J} M_n(\cl S\otimes_{\tau_j}\cl T)^+$. Since
$P_n\subseteq M_n(\cl S\otimes_{\tau_{j_0}}\cl T)^+$ for each
$j_0\in J$, it follows that $P_n \cap (-P_n) = \{0\}.$ It is trivial
to check that the family $\{P_n\}_{n=1}^{\infty}$ is compatible and
that it satisfies $M_n(\cl S)^+ \otimes M_m(\cl T)^+ \subseteq
P_{mn}.$  Hence $(P_n - P_n) + i(P_n - P_n) = M_n(\cl S \otimes \cl
T).$ Thus $\{P_n \}_{n=1}^{\infty}$ is a matrix ordering on $\cl S
\otimes \cl T$. We shall denote this matrix-ordered space by $\cl S
\otimes_{\tau} \cl T$.

Since $M_n(\cl S \otimes_{\tau_j} \cl T)^+\subseteq M_n(\cl S \otimes_{\min} \cl T)^+$ for every $j \in J,$ it follows that
$P_n\subseteq M_n(\cl S \otimes_{\min} \cl T)^+$, $n\in \bb{N}$.
Since $1 \otimes 1$ is a matrix order unit for $\cl S \otimes_{\min} \cl T,$
it follows that $1 \otimes 1$ is a matrix order unit for $\cl S \otimes_{\tau}
\cl T.$  Also, since $1 \otimes 1$ is Archimedian for each $\cl S
\otimes_{\tau_j} \cl T$, it follows that $1 \otimes 1$ is
Archimedian for $\cl S \otimes_{\tau} \cl T$.  Hence
$\cl S \otimes_{\tau} \cl T$ is an operator system, that is, Property~(T1) holds. The fact
that Property~(T2) holds follows from the fact that $M_n(\cl S \otimes_{\max} \cl T)^+
\subseteq M_n(\cl S \otimes_{\tau} \cl T)^+$.  Property~(T3) holds
because it holds $\min$ and
$M_n(\cl S \otimes_{\tau} \cl T)^+ \subseteq M_n(\cl S
\otimes_{\min} \cl T)^+$.

Finally, if every $\tau_j$ is functorial and $\phi_i: \cl S_i \to \cl
T_i$ for $i=1,2$ are unital completely
positive maps,
then $\phi_1 \otimes \phi_2: \cl S_1 \otimes_{\tau_j} \cl T_1 \to \cl
S_2 \otimes_{\tau_j} \cl T_2$ is a unital completely
positive map for every $j \in J$. Since
the positive cones for $\cl S_1 \otimes_{\tau} \cl T_1$ are smaller than the positive cones for $\cl S_1 \otimes_{\tau_j} \cl T_1$,
we have that $\phi_1 \otimes \phi_2: \cl S_1 \otimes_{\tau} \cl T_1
\to \cl S_2 \otimes_{\tau_j} \cl T_2$ is a unital completely
positive map for every $j \in J.$
From this it follows that $\phi_1 \otimes \phi_2: \cl S_1
\otimes_{\tau} \cl T_2 \to \cl S_2 \otimes_{\tau} \cl T_2$ is a unital completely
positive map, and the functoriality of $\tau$ follows.
\endproof

Motivated by the previous section, we introduce a general way to
induce operator system structures from inclusions. Let $\alpha$ be an operator system tensor product.
If $\cl S_i$ and $\cl T_i$, $i = 1,2$, are operator systems with $\cl S_1\subseteq \cl S_2$ and $\cl T_1\subseteq \cl T_2$, let $\{C_n\}_{n=1}^{\infty}$ be the matrix ordering on $\cl S_1\otimes\cl T_1$ given by
$$C_n = M_n(\cl S_2\otimes_{\alpha}\cl T_2)^+ \cap M_n(\cl S_1\otimes\cl T_1), \qquad \text{ $n\in \bb{N}$}.$$
We call $\{C_n\}_{n=1}^{\infty}$ the operator system structure on $\cl
S_1\otimes\cl T_1$ {\bf induced by $\alpha$} and the pair $(\cl S_2,\cl T_2)$. We note that
this is not an operator system tensor product in the sense of
definition given in Section~\ref{s-gs}; it is defined ``locally'' for every quadruple of operator systems $\cl S_1\subseteq
\cl S_2$ and $\cl T_1\subseteq\cl T_2$.

A tensor product $\alpha$ on the category of operator systems is
called {\bf left injective} if for all operator systems $\cl S_1$,
$\cl S_2$, and $\cl T$ with $\cl S_1\subseteq \cl S_2$, the
inclusion of $\cl S_1\otimes_{\alpha}\cl T$ into $\cl
S_2\otimes_{\alpha}\cl T$ is a complete order isomorphism. Equivalently,
$\alpha$ is left injective if the operator system structure of
$\cl S_1\otimes_{\alpha}\cl T$
coincides with the one induced by $\alpha$ and $(\cl
S_2,\cl T)$ for all operator systems $\cl S_2$ with $\cl
S_1\subseteq\cl S_2$, and all operator systems $\cl T$. We define a {\bf right injective} operator
system tensor product similarly.
 An operator system
tensor product is {\bf injective} if it is both left and right
injective or, equivalently, if the inclusion of $\cl S_1
\otimes_{\alpha} \cl T_1$ into $\cl S_2 \otimes_{\alpha} \cl T_2$ is a
complete order injection whenever $\cl S_1 \subseteq \cl S_2$ and $\cl
T_1 \subseteq \cl T_2.$  For example, $\min$ is an injective tensor
product.

Given an operator system $\cl S$ we let $I(\cl S)$ denote its
injective envelope. There is a precise sense in which $I(\cl S)$ is
the ``smallest'' injective operator system that contains $\cl S.$ See
\cite[Chapter 15]{Pa} for a detailed development of this concept.

\begin{definition}\label{d_lre}
Let $\cl S$ and $\cl T$ be operator systems. We let $\cl
S\otimes_{\linj}\cl T$ (respectively, $\cl S\otimes_{\rinj}\cl T$) be the operator system with underlying space
$\cl S\otimes\cl T$ whose matrix ordering is
induced by the inclusion $\cl S\otimes\cl T\subseteq I(\cl
S)\otimes_{\max}\cl T$ (respectively, $\cl S\otimes\cl T\subseteq \cl
S\otimes_{\max} I(\cl T)$).

Likewise, we let $\cl S\otimes_{\inj}\cl T$ be the operator system with underlying space
$\cl S\otimes\cl T$ whose matrix ordering is
induced by the inclusion $\cl S\otimes\cl T\subseteq I(\cl
S)\otimes_{\max} I(\cl T)$.
\end{definition}

\begin{theorem}\label{th_linjf}
The mappings $\linj : \cl O \times \cl O \to \cl O$ $\rinj: \cl O
\times \cl O \to \cl O$ and $\inj: \cl O \times \cl O \to \cl O$ sending the pair
$(\cl S,\cl T)$ to the operator system $\cl S\otimes_{\linj}\cl T,$
$\cl S \otimes_{\rinj} \cl T$ and $\cl S \otimes_{\inj} \cl T$
are functorial operator system tensor products.
\end{theorem}
\begin{proof} We only prove the left injective case, the other proofs are
  similar.

Properties (T1) and (T2) are immediate from the definition of $\linj$ and
the fact that $\max$ is an operator system tensor product.
Let $\cl S$ and $\cl T$ be operator systems.
Suppose that $\phi\in S_n(\cl S)$ and $\psi\in S_m(\cl T)$, and let $\tilde{\phi}\in S_n(I(\cl S))$
be an extension of $\phi$. Since $\max$ is an operator system tensor product, by (T3) we have that
$\tilde{\phi}\otimes\psi : I(\cl S)\otimes_{\max} \cl T\rightarrow M_{mn}$ is completely positive.
Restricting to $\cl S\otimes_{\linj}\cl T$, we obtain that $\phi\otimes\psi : \cl S\otimes_{\linj}\cl T\rightarrow M_{mn}$
is completely positive. Thus, $\linj$ possesses Property (T3).

Now let $\cl S_i$ and $\cl T_i$ be operator systems, $i = 1,2$, and
$\phi \in {\rm CP}(\cl S_1,\cl S_2)$, $\psi \in {\rm CP}(\cl T_1,\cl T_2)$.
Let $\tilde{\phi} : I(\cl S_1)\rightarrow I(\cl S_2)$ be a completely positive extension of
$\phi$. By the functoriality of $\max$, we have that
$\tilde{\phi}\otimes\psi : I(\cl S_1)\otimes_{\max}\cl T\rightarrow I(\cl S_2)\otimes_{\max} \cl T$
is completely positive. Restricting to $\cl S_1\otimes_{\linj}\cl T$, we obtain that
$\phi\otimes\psi : \cl S_1\otimes_{\linj}\cl T_1\rightarrow \cl S_2\otimes_{\linj}\cl T_2$
is completely positive.
\end{proof}

\begin{lemma}\label{l_slin}
Let $\cl S$, $\cl S_1$, and $\cl T$ be operator systems with $\cl
S\subseteq\cl S_1$, and let $\tau$ be the operator
system structure induced by the inclusion $\cl S\otimes\cl T\subseteq \cl
S_1\otimes_{\max}\cl T$. Then $\cl S\otimes_{\tau} \cl T$ is greater
than $\cl S\otimes_{\linj}\cl T$.
\end{lemma}
\begin{proof}
Let $\phi : \cl S_1\rightarrow I(\cl S)$ be a unital completely
positive map extending the inclusion
$\iota : \cl S\rightarrow I(\cl S)$. By the functoriality of the maximal operator system
tensor product, we have that $\phi\otimes\id : \cl S_1\otimes_{\max}\cl T\rightarrow I(\cl S)\otimes_{\max} \cl T$
is completely positive. Since $\phi\otimes\id$ coincides on $\cl S\otimes\cl T$ with the identity map,
the conclusion follows.
\end{proof}

We now show the role that these tensor products play within the family
of all operator system tensors.

\begin{theorem}\label{th_smlinj}
The operator system tensor product $\linj$ is left injective.
Moreover, if $\alpha: \cl O \times\cl O\rightarrow\cl O$ is a left
injective functorial operator system tensor product then $\linj$ is
greater than $\alpha$.
 Similarly, $\rinj$ is the largest right
injective and $\inj$ is the largest injective functorial operator
system tensor products.
\end{theorem}
\begin{proof} We only prove the first statement.

Suppose that $\cl S\subseteq\cl S_1$. Let $\cl
S\otimes_{\tau} \cl T$ denote the operator system induced by the inclusion
$\cl S\otimes\cl T\subseteq I(\cl S_1)\otimes_{\max}\cl T$. By Lemma~\ref{l_slin},
$\cl S\otimes_{\tau} \cl T$ is greater than $\cl S\otimes_{\linj} \cl T$.
On the other hand,
the inclusion $\cl S\subseteq I(\cl S_1)$ gives rise to a unital completely
positive map
$\phi : I(\cl S)\rightarrow I(\cl S_1)$. By functoriality, the map
$\phi\otimes\id : I(\cl S)\otimes_{\max}\cl T\rightarrow I(\cl
S_1)\otimes_{\max}\cl T$ is completely positive. Restricting to the
subspace $\cl S\otimes\cl T$ implies that the corresponding map
$\phi\otimes\id : \cl S\otimes_{\linj}\cl T\rightarrow \cl
S\otimes_{\tau}\cl T$ is completely positive.  Since this map coincides with the identity map, we have that $\cl
S\otimes_{\linj}\cl T$ is greater than $\cl S\otimes_{\tau}\cl T$,
and hence $\cl S \otimes_{\tau} \cl T = \cl S\otimes_{\linj}\cl T$.
Thus the inclusion $\cl S\otimes_{\linj}\cl T\subseteq \cl
S_1\otimes_{\linj}\cl T$ is completely isometric. It is thus shown
that $\linj$ is injective.

Suppose now that $\alpha : \cl O\times \cl O\rightarrow\cl O$ is a
left injective operator system tensor product. If $\cl S$ and $\cl
T$ are operator systems, then $\cl S\otimes_{\alpha} \cl T \subseteq
I(\cl S)\otimes_{\alpha} \cl T$ completely order
isomorphically. By the maximality property of $\max$, we have that
the identity map $\id\otimes\id : I(\cl S)\otimes_{\max}\cl
T\rightarrow I(\cl S)\otimes_{\alpha}\cl T$ is completely positive.
Hence its restriction to $\cl S\otimes\cl T$ maps the positive cones
of $\cl S\otimes_{\linj}\cl T$ into those of $\cl
S\otimes_{\alpha}\cl T$. Thus $\linj$ is greater than $\alpha$.
\end{proof}

We summarize the order relations between the particular tensor products studied in this paper:
$$\min \leq \inj\leq \linj,\rinj \leq \comm \leq \max.$$

Since $\linj$ and $\rinj$ play central roles in the family of all
tensor products, it is interesting to know if their
relationship to important properties of C*-algebras. The following
results provide partial answers to these questions.

\begin{proposition}\label{th_wep}
Let $A$ be a unital C*-algebra. The following are equivalent:
\begin{enumerate}
\item[(i)] $A$ possesses the weak expectation property (WEP);
\item[(ii)] $A\otimes_{\linj} B = A\otimes_{\max} B$ for every C*-algebra
$B$.
\end{enumerate}
\end{proposition}
\begin{proof}
(i)$\Rightarrow$(ii) By Lance's characterization of WEP (see \cite{La}),
 the inclusion of $A\otimes_{\max} B$ into
$I(A)\otimes_{\max} B$ is a complete order isomorphism onto its range. However, $A\otimes_{\linj} B$ is by
definition obtained by restricting the matrix order structure of
$I(A)\otimes_{\max} B$ to $A\otimes B$. It follows that
$A\otimes_{\max} B = A\otimes_{\linj} B$.

(ii)$\Rightarrow$(i) Suppose that $A_1$ and $B$ are C*-algebras such
that $A\subseteq A_1$. By Lemma \ref{l_slin}, the matrix ordering on
$A\otimes B$ induced by its inclusion in $A_1\otimes_{\max} B$ is
(set-theoretically) contained in that of $A\otimes_{\linj} B =
A\otimes_{\max} B$. However, it is trivial that the matrix ordering
of $A\otimes_{\max} B$ is contained in the former matrix ordering
since $A\otimes_{\max} B$ is the largest matrix ordering on
$A\otimes B$. Thus, $A\otimes_{\max} B \subseteq A_1\otimes_{\max}
B$ (as C*-algebras). It follows from \cite{La} that $A$ has WEP.
\end{proof}

Proposition \ref{th_wep} shows that WEP can be thought of as a
nuclearity property with respect to $\linj$, which is an operator
system structure on the tensor products bigger than the minimal one.
The next observation characterizes nuclearity in terms of the right
injective tensor product $\rinj$.

\begin{proposition}\label{p_nuc}
Let $A$ be a unital C*-algebra. The following are equivalent:
\begin{enumerate}
\item[(i)] $A$ is nuclear,

\item[(ii)] $A\otimes_{\rinj} B = A\otimes_{\max}B$, for every unital C*-algebra $B$.
\end{enumerate}
\end{proposition}
\proof
(i)$\Rightarrow$(ii) If $A$ is nuclear then
$A\otimes_{\min} B = A\otimes_{\max} B$ sits completely order isomorphically in
$A\otimes_{\min} I(B) = A\otimes_{\max} I(B)$, and hence $A\otimes_{\max} B = A\otimes_{\rinj} B$.

(ii)$\Rightarrow$(i) Let $B$ and $B_1$ be unital C*-algebras with $B\subseteq B_1$.
Let $\phi : B_1\rightarrow I(B)$ be a completely positive extension of the inclusion $B\rightarrow I(B)$.
Suppose that $u\in M_n(A\otimes B) \cap M_n(A\otimes_{\max} B_1)^+$.
Using the identifications
$M_n(A\otimes B) \equiv M_n(A)\otimes B$ and
$M_n(A\otimes_{\max} B_1) \equiv M_n(A)\otimes_{\max} B_1$ and the functoriality
of the maximal tensor product, we have that
$$(\id\mbox{}_{M_n(A)}\otimes \phi)(u) \in (M_n(A)\otimes_{\max} I(B))^+ \equiv M_n(A\otimes_{\max} I(B))^+.$$
Since $u\in M_n(A\otimes B)$ and $\phi$ coincides with the identity mapping on $B$,
we have that $u\in M_n(A\otimes_{\max} I(B))^+$. By assumption,
$u\in M_n(A\otimes_{\max} B)^+$. We thus showed that
the inclusion $A\otimes_{\max} B\rightarrow A\otimes_{\max} B_1$ is a complete order
isomorphism onto its range. It follows from \cite[Theorem A]{La} that $A$ is nuclear.
\endproof

Proposition \ref{p_nuc} allows one to establish the nuclearity of a
C*-algebra by comparing the maximal tensor product with $\rinj$,
which is a priori bigger than the minimal tensor product.

Propositions \ref{th_wep} and \ref{p_nuc} have the following consequence.

\begin{corollary}\label{c_assym}
The tensor product $\linj$ is not symmetric.
\end{corollary}
\proof By \cite{La}, there exists a C*-algebra $A$ which is not
nuclear and possesses the weak expectation property. By Propositions
\ref{th_wep} and \ref{p_nuc}, there exists a unital C*-algebra $B$
such that $A\otimes_{\rinj} B \neq A\otimes_{\max} B =
A\otimes_{\linj} B.$

Suppose that the map $\theta : A\otimes B\rightarrow B\otimes A$
given by $\theta(x\otimes y) = y \otimes x$ was a complete order
isomorphism of $A\otimes_{\linj} B$ onto $B\otimes_{\linj} A$. Since
$A$ has WEP, Proposition \ref{th_wep} implies that $\theta :
A\otimes_{\max} B \rightarrow I(B)\otimes_{\max} A$ is a complete
order isomorphism onto its range. Since $\max$ is symmetric, the
restriction of the mapping $\theta^{-1} : I(B)\otimes_{\max}
A\rightarrow A\otimes_{\max} I(B)$ to $B\otimes A$ is a complete
order isomorphism onto its range. It follows that the inclusion
$A\otimes_{\max} B\rightarrow A\otimes_{\max} I(B)$ is a complete
order isomorphism onto its range, and hence $A\otimes_{\max} B =
A\otimes_{\rinj} B$, a contradiction with the choice of $B$.
\endproof



\begin{remark}
{\rm Arguments similar to those given above
show that if $X$ and $Y$ are operator spaces,
then the inclusion $X\otimes Y\subseteq I(X)\hat{\otimes}I(Y)$
induces an operator space tensor product $X\otimes_{\hat{\inj}}Y$
that is the largest injective tensor product in the operator space category.
We claim that the operator space structure on $X\otimes_{\hat{\inj}}Y$
is distinct from the one on $X\otimes^{\inj} Y$ (recall that
$X\otimes^{\inj} Y$ arises from the embedding $X\otimes Y\subseteq
\cl S_X\otimes_{\inj} \cl S_Y$ --- or, equivalently, from the embedding
$X\otimes Y\subseteq I(\cl S_X)\otimes_{\max} I(\cl S_Y)$). To see
this, let $X = Y = M_{m,n}$. Then $I(\cl S_X) = M_{m+n,m+n}$ and by
the nuclearity of $M_{m+n,m+n}$ we have that $I(\cl
S_X)\otimes_{\max} I(\cl S_Y) = M_{(m+n)^2,(m+n)^2}$. However,
$M_{m,n}\hat{\otimes}M_{m,n}$ is distinct from $M_{m^2,n^2}$. Hence,
$X\otimes_{\hat{\inj}}Y\neq X\otimes^{\inj} Y$ in this case. As a
corollary we obtain the following.}
\end{remark}

\begin{corollary}
There exists a functorial injective operator space tensor product
that is not induced by a functorial injective operator system
tensor product.
\end{corollary}

In our last proposition, we characterize the norm $\|\cdot\|_{\inj}$
induced by the operator system structure $\inj$ introduced in
Definition \ref{d_lre}.

\begin{proposition}\label{p_chareno}
Let $A$ and $B$ be unital C*-algebras and $u\in A\otimes B$. Then
\begin{align*}
\|u\|_{\inj} = \inf\{\|u\|_{A_1\otimes_{\max} B_1} : \text{$A_1$ and }&\text{ $B_1$ are C*-algebras} \\
&\text{with $1_A\in A\subseteq A_1$ and $1_B\in B\subseteq B_1$ } \}.
\end{align*}
\end{proposition}
\proof Fix $u\in A\otimes B$ and denote the quantity on the right
hand side by $\delta$. By the definition of $\inj$ and $\delta$, we
have that $\delta \leq \|u\|_{\inj}$.

Let $A_1$ and $B_1$ be C*-algebras with $1_A\in A\subseteq A_1$ and
$1_B\in B\subseteq B_1$. Let $\phi : A_1\rightarrow I(A)$ and $\psi :
B_1\rightarrow I(B)$ be completely positive extensions of the
inclusion maps $A\rightarrow I(A)$ and $B\rightarrow I(B)$,
respectively. By functoriality, $\phi\otimes\psi$ is a unital
completely positive, and hence completely contractive, map from
$A_1\otimes_{\max}B_1$ into $I(A)\otimes_{\max} I(B)$. It follows
that
$$\|u\|_{\inj} = \|(\phi\otimes\psi)(u)\|_{I(A)\otimes_{\max}I(B)}
\leq \|u\|_{A_1\otimes_{\max} B_1},$$ and hence $\delta =
\|u\|_{\inj}$.
\endproof

\begin{remark}{\rm  Pisier \cite[p.~350]{Pi} defines a tensor product $\otimes_M$
    on operator spaces $X \subseteq B(H)$ and $Y \subseteq B(K)$ by
    identifying $X \otimes_M Y$ with the subspace $X \otimes Y \subseteq
    B(H) \otimes_{\max} B(K),$ and argues that this tensor product is
    independent of the particular completely isometric inclusions of $X$ and $Y$ into
    $B(H)$ spaces. It is not
    difficult to see that this tensor product is identical with our
    tensor product $\otimes^e$. We make this precise in the
    following.}
\end{remark}

Recall that every operator system is also an operator space. Thus,
we may form the operator system $\cl S \otimes_e \cl T$ and the
operator space $\cl S \otimes^e \cl T.$

\begin{proposition} Let $X$ and $Y$ be operator spaces and let $\cl S$
  and $\cl T$ be operator systems. Then $X \otimes^e Y = X \otimes_M
  Y$ and $\cl S \otimes_e \cl T = \cl S \otimes^e \cl T,$ completely
  isometrically.
\end{proposition}
\begin{proof} First let $I_H \in A \subseteq \cl B(H)$ and $I_K \in  \cl B
  \subseteq  \cl B(K)$ be unital,
  injective C*-subalgebras. Then there exists unital completely positive
  projections $\phi : \cl B(H) \to A$ and $\psi : \cl B(K) \to B.$ This implies
  that the map $\phi \otimes \psi : \cl B(H) \otimes_{\max}  \cl B(K) \to A
  \otimes_{\max} B,$ is a unital completely positive map. Hence it
  follows that the operator subsystem $A \otimes B \subseteq \cl B(H)
  \otimes_{\max}  \cl B(K)$ is completely order isomorphic to $A
  \otimes_{\max} B.$

Thus if we are given operator spaces $X$ and $Y$ and we embed
$I(S_X) \subseteq  \cl B(H)$ and $I(S_Y) \subseteq  \cl B(K),$ then
the subspaces $X \otimes^e Y \subseteq I(S_X) \otimes_{\max} I(S_Y)$
and $X \otimes_M Y \subseteq  \cl B(H) \otimes_{\max}  \cl B(K),$
will be completely isometric.

If $\cl S \subseteq I(\cl S) \subseteq  \cl B(H)$ and $\cl T
\subseteq I(\cl T) \subseteq  \cl B(K)$ are operator systems, then
the previous paragraph shows that $\cl S \otimes^e \cl T = \cl S
\otimes_M \cl T$ completely isometrically. But $\cl S \otimes_M \cl
T$ can be completely isometrically identified with the subspace $\cl
S \otimes \cl T \subseteq  \cl B(H) \otimes_{\max}  \cl B(K),$ and
we also have the completely isometric identification $I(\cl S)
\otimes_{\max} I(\cl T) \subseteq  \cl B(H) \otimes_{\max}  \cl
B(K).$ Hence we have that $\cl S \otimes_M \cl T \subseteq I(\cl S)
\otimes_{\max} I(\cl T)$ is a completely isometric inclusion and so
$\cl S \otimes_M \cl T = \cl S \otimes_e \cl T.$
\end{proof}

In contrast, recall that even if $A$ and $B$ are unital C*-algebras, then
$A \otimes^{\max} B = A \hat{\otimes} B$, which is not completely
isometrically equal to $A \otimes_{\max} B.$

\medskip

\noindent {\bf Acknowledgements. } We would like to thank K.H. Han
for his careful reading of the manuscript which lead to many
improvements in the presentation, and E.G. Effros and G. Pisier for
helpful comments and references.

\end{document}